\documentclass[a4paper, 11pt]{article}

\usepackage[utf8]{inputenc}   

\usepackage[centertags]{amsmath}
\usepackage{amssymb,amsfonts,amsthm}
\usepackage{indentfirst}
\usepackage{array}
\usepackage{float}
\usepackage[usenames,dvipsnames,svgnames,table]{xcolor}
\usepackage{hyperref}
\usepackage{enumerate}
\usepackage{graphicx}
\usepackage{epstopdf}
\usepackage{graphics}
\usepackage{subfigure}
\usepackage[margin=2.5cm,footskip=1.5cm,includefoot]{geometry} 
\usepackage{listings}
\usepackage{setspace}
\usepackage{multirow}
\usepackage{caption}
\usepackage[backend=biber,style=numeric,sorting=nyt,backref]{biblatex}
\usepackage{amsmath,enumitem,amsthm,amssymb,graphicx,mathtools,tikz,yfonts,tikz-cd, subfiles, datetime, adjustbox, extarrows,old-arrows, algorithm2e}
    \usetikzlibrary{decorations.pathreplacing,angles,quotes,positioning}

\usepackage[nottoc]{tocbibind}

\DeclareFieldFormat{postnote}{#1}
\DeclareFieldFormat{multipostnote}{#1}

\usetikzlibrary{positioning,lindenmayersystems}

\usepackage{graphicx}
\graphicspath{ {./images/} }
\usepackage[rightcaption]{sidecap}
\usepackage{wrapfig}

\usepackage[capitalize,nameinlink,noabbrev]{cleveref}

\definecolor{darkgreen}{HTML}{00563F}
\definecolor{pinkish}{HTML}{E80072}
\hypersetup{
    colorlinks  = true,
    linkcolor   = blue,
    filecolor   = magenta,      
    urlcolor    = cyan,
    citecolor   = pinkish,
}

\theoremstyle{plain}

\newtheorem{theorem}{Theorem}[section]

\newtheorem{proposition}[theorem]{Proposition}
\AddToHook{env/proposition/begin}{\crefalias{theorem}{proposition}}

\newtheorem{conjecture}[theorem]{Conjecture}
\AddToHook{env/conjecture/begin}{\crefalias{theorem}{conjecture}}

\newtheorem{corollary}[theorem]{Corollary}
\AddToHook{env/corollary/begin}{\crefalias{theorem}{corollary}}

\newtheorem{lemma}[theorem]{Lemma}
\AddToHook{env/lemma/begin}{\crefalias{theorem}{lemma}}

\newtheorem{definition}[theorem]{Definition}
\AddToHook{env/definition/begin}{\crefalias{theorem}{definition}}

\newtheorem*{convention}{Convention}

\DeclareMathOperator{\id}{id}

\renewcommand{\div}{\mathrm{div}}

\newcommand{\R}{\mathbb R}
\newcommand{\C}{\mathbb C}
\newcommand{\Z}{\mathbb Z}

\newcommand{\Q}{\mathbb Q}
\renewcommand{\P}{\mathbb P}

\newcommand{\conv}{\mathrm{Conv}}
\newcommand{\cone}{\mathrm{Cone}}
\newcommand{\vol}{\mathrm{vol}}
\newcommand{\nvol}{\mathrm{n\text{-}vol}}
\DeclareFontFamily{U}{mathx}{}
\DeclareFontShape{U}{mathx}{m}{n}{<-> mathx10}{}
\DeclareSymbolFont{mathx}{U}{mathx}{m}{n}
\DeclareMathAccent{\widecheck}{0}{mathx}{"71}
\makeatletter
\newcommand{\tpitchfork}{%
  \vbox{
    \baselineskip\z@skip
    \lineskip-.52ex
    \lineskiplimit\maxdimen
    \m@th
    \ialign{##\crcr\hidewidth\smash{$-$}\hidewidth\crcr$\pitchfork$\crcr}
  }%
}
\makeatother

\newcommand{\parens}[1]{{\left(#1\right)}}
\newcommand{\bracket}[1]{{\left[#1\right]}}
\newcommand{\curly}[1]{{\left\{#1\right\}}}
\newcommand{\Angle}[1]{{\left\langle#1\right\rangle}}
\newcommand{\pipe}[1]{{\left|#1\right|}}

\begin{document}
    \title{Smooth Calabi-Yau varieties with large index and Betti numbers}
    \author{Jas Singh\\
            \footnotesize{University of California, Los Angeles}\\
            \footnotesize{\hyperlink{mailto:jas@math.ucla.edu}{jas@math.ucla.edu}}}
    \date{  }

    \maketitle
    
	\begin{abstract}
        \footnotesize
        {
            We construct smooth Calabi-Yau varieties in every dimension with doubly exponentially growing index, which we conjecture to be maximal in every dimension.
            We also construct smooth Calabi-Yau varieties with extreme topological invariants; namely, their Euler characteristics and the sums of their Betti numbers grow doubly exponentially.
            These are conjecturally extremal in every dimension.
            The varieties we construct are known in small dimensions, but we believe them to be new in general.
            This work builds off of the singular Calabi-Yau varieties found by Esser, Totaro, and Wang in \cite{largeindex}.
        }
	\end{abstract}

	\section{Introduction}
        \label{section/intro}
        We work over the complex numbers.
A  normal, projective variety $X$ is called \emph{Calabi-Yau} if its canonical divisor satisfies $K_X \sim_{\Q} 0$.
In other words, $K_X$ is a torsion element in the divisor class group of $X$.
Note that other authors insist that Calabi-Yau varieties are simply connected, or that $H^i(X, \mathcal O_X) = 0$ for $0 < i < \dim X$.
We do not insist on these additional constraints.

The index of a Calabi-Yau variety $X$ is the torsion order of $K_X$, i.e. the smallest positive $m$ so that $m K_X \sim 0$.
We are interested in how large the index of smooth, projective Calabi-Yau $n$-folds can be.

Let $s_n$ be Sylvester's sequence, defined by $s_0 = 2$, $s_{n + 1} = s_0 \cdots s_n + 1$, which grows doubly exponentially.
Our first main result is the following.
\begin{theorem} \label{intro/largeindex}
    For $n \geq 1$ there exists a smooth, projective Calabi-Yau $n$-fold $V^{(n)}$ with index $(s_{n - 1} - 1)(2 s_{n - 1} - 3)$.
\end{theorem}
\begin{center}
    \begin{tabular}{c|c}
        $n$ & $\mathrm{index}(V^{(n)})$\\
        \hline
        1 & 1\\
        2 & 6\\
        3 & 66\\
        4 & 3486\\
        5 & 6521466\\
        6 & 21300104111286\\
        7 & 226847426110811738551148466
    \end{tabular}
\end{center}

In small dimensions, these examples are known to have the largest possible index, as proven for 3-folds in \cite[Corollary 5]{k3}.
For instance, in dimension 1, any elliptic curve will do.
In dimension 2, $V^{(2)} = (E_1 \times E_2)/\mu_6$, where $E_1$ is the elliptic curve whose automorphism group is the cyclic group $\mu_6$ (which is the largest possible), and where $E_2$ is any elliptic curve acted on by $\mu_6$ via translation by a $6$-torsion point.
In dimension 3, $V^{(3)} = (S \times E)/\mu_{66}$, where $S$ is the K3 surface with a non-symplectic automorphism of order $66$ (which is the largest possible), and $E$ is an elliptic curve acted on via translation by a $66$-torsion point.

A major conjecture due to Shokurov is that the indices of Calabi-Yau varieties are bounded in each dimension.
More precisely, Xu formulated the following conjecture in \cite[Conjecture 1.8.3]{indexconjecture}.
\begin{conjecture}[Index Conjecture]
    For a positive integer $n$ and a subset $I \subseteq [0, 1]$ satisfying the descending chain condition, there is a positive integer $c(n, I)$ so that for any klt Calabi-Yau pair $(X, D)$ of dimension $n$ with the coefficients of $D$ lying  in $I$, we have $c(n, I)(K_X + D) \sim 0$.
\end{conjecture}
A pair $(X, D)$ is Calabi-Yau if $K_X + D \sim_{\Q} 0$, where $X$ is a normal, projective variety and $D$ is an effective $\Q$-divisor on $X$.
The index of the pair is defined to be the torsion order of $K_X + D$.
The index conjecture thereby asserts a uniform upper bound to the indices of klt Calabi-Yau pairs with fixed dimension $n$ and (dcc) set of allowable coefficients $I$ for the boundary divisor.

The conjecture was proven for surfaces in \cite[Corollary 1.11]{ps09}.
In dimension 3, the index conjecture has been proven for terminal 3-folds in \cite[Theorem 3.2]{kawamata86} and \cite[Corollary]{morrison86}, klt 3-folds in \cite[Corollary 1.7]{jiang21}, lc pairs \cite[Theorem 1.8.7]{indexconjecture}, and slc pairs in \cite[Corollary 1.6]{jl21}.
In dimenson 4, \cite[Corollary 1.7]{jl21} shows the index conjecture for lc pairs which are not klt and \cite[Theorem 1.8.7]{indexconjecture} shows it for lc pairs with nonzero boundary which are klt.
Most recently, the index conjecture has been proven for rationally connected varieties by Birkar in \cite[Theorem 1.6]{birkar2025}.
Note that the Calabi-Yau pairs with large index defined by Esser, Totaro, and Wang in \cite{largeindex} are rationally connected, but the $V^{(n)}$ of \cref{intro/largeindex} are not.

This conjecture remains open, so it remains unclear if there is any uniform upper bound to the indices of Calabi-Yau varieties of a fixed dimension, even with a fixed set of allowable coefficients.
Nevertheless, Esser, Totaro, and Wang conjectured in \cite{largeindex} that $(s_{n - 1} - 1)(2 s_{n - 1} - 3)$ is precisely that upper bound for terminal Calabi-Yau $n$-folds, and they found explicit examples achieving that value.
We formulate therefore a similar conjecture in the smooth case, which follows from theirs.
\label{intro/smoothindex}
\begin{conjecture}
    The largest possible index of any smooth, projective Calabi-Yau $n$-fold is $(s_{n - 1} - 1) (2 s_{n - 1} - 3)$.
\end{conjecture}

Let $I_{\mathrm{sm}}(n)$, $I_{\mathrm{term}}(n)$ be the sets of indices of smooth, terminal Calabi-Yau $n$-folds, respectively.
Additionally, let $I_{\mathrm{klt \; st}}(n)$ be the set of indices of klt Calabi-Yau pairs with \emph{standard coefficients} of dimension $n$.
(We say that a pair $(X, D)$ has standard coefficients if the coefficients of $D$ are all of the form $1 - \frac{1}{d}$ for $d$ a positive integer.
So, in the notation of the index conjecture, we are taking $I = \{1 - \frac{1}{d} \mid d \in \Z^{>0}\}$.)

In \cite{largeindex}, Esser, Totaro, and Wang constructed a klt Calabi-Yau pair with standard coefficients of dimension $n - 1$ with large index, and then used it to construct a terminal Calabi-Yau variety of dimension $n$ with the same index.
In fact, it's true that $I_{\mathrm{klt \; st}}(n - 1) \subseteq I_{\mathrm{term}}(n)$, and Esser, Totaro, and Wang's conjecture can be rephrased as saying that both sets have the same maximum: $(s_{n - 1} - 1)(2 s_{n - 1} - 3)$.
Masamura conjectured further in \cite{smoothindices} that $I_{\mathrm{klt\;st}}(n - 1) = I_{\mathrm{term}}(n)$ and proved that $I_{\mathrm{sm}}(n) \subseteq I_{\mathrm{klt\;st}}(n - 1) \subseteq I_{\mathrm{term}}(n)$.
Given Esser, Totaro, and Wang's conjectural upper bound on the indices of terminal Calabi-Yau $n$-folds, our \cref{intro/largeindex} would show that $I_{\mathrm{sm}}(n)$ and $I_{\mathrm{term}}(n)$ share the same maximum, which is consistent with Masamura's conjecture.

The construction of the $V^{(n)}$ in \cref{intro/largeindex} will also lead to the definition of smooth Calabi-Yau varieties with extreme topological invariants.
\begin{theorem}
\label{intro/largebetti}
    For every $n \geq 1$, there exists a smooth, projective Calabi-Yau $n$-fold $W^{(n)}$ with the following properties:
    \begin{itemize}
        \item The sum of the Betti numbers of $W^{(n)}$ is $2 (s_0 - 1) \cdots (s_n - 1)$.
        \item When $n$ is even, the Euler characteristic of $W^{(n)}$ is $2 (s_0 - 1) \cdots (s_n - 1)$.
        \item When $n$ is odd, the Euler characteristic of $W^{(n)}$ is $-(s_0 - 1) \cdots (s_{n - 1} - 1) (2 s_n - 6)$.
    \end{itemize}
\end{theorem}
\begin{center}
    \begin{tabular}{c|c|c}
        $n$ & $\sum_{i = 0}^{2n} b_i(W^{(n)})$ & $\chi(W^{(n)})$\\
        \hline
        1 & 4 & 0\\
        2 & 24 & 24\\
        3 & 1008 & -960\\
        4 & 1820448 & 1820448\\
        5 & 5940926462016 & -5940922821120\\
        6 & 63271205161020798539584896 & 63271205161020798539584896
    \end{tabular}
\end{center}

We now elaborate on our strategy to construct the varieties $V^{(n)}$ and $W^{(n)}$ of our main theorems.

Let's start with $V^{(n)}$.
Let $m = (s_{n - 1} - 1)(2 s_{n - 1} - 3)$.
Esser, Totaro, and Wang constructed their klt Calabi-Yau pair of large index by defining a hypersurface $X$ in an $n$-dimensional weighted projective space $\P$ which is stabilized by a cyclic subgroup $\mu_{m}$ of toric automorphisms of $\P$, and taking the quotient pair $X/\mu_m$.
To construct their terminal example from this, take $Z \longrightarrow X$ a $\mu_{m}$-equivariant terminalization, and let $E$ be an elliptic curve with a fixed $m$-torsion point.
Then $\mu_m$ acts on $E$ via translation by this point and $(Z \times E)/\mu_{m}$ is a terminal Calabi-Yau $n$-fold with index $m$.
If we can instead find a $\mu_{m}$-equivariant, crepant resolution of singularities $\widetilde{X} \longrightarrow X$, the same procedure will yield a smooth Calabi-Yau $n$-fold of index $m$:
\[
    V^{(n)} = \frac{\widetilde{X} \times E}{\mu_m}.
\]

As for $W^{(n)}$, Esser, Totaro, and Wang prove that the hypersurface $X \subseteq \P$ mentioned in the construction of $V^{(n)}$ has the values in \cref{intro/largebetti} as its \emph{orbifold} Betti numbers, as defined by Chen and Ruan in \cite{orbifold}.
In \cite{orbifoldcrepant}, Yasuda proved that for a projective variety with Gorenstein quotient singularities (such as $X$), its orbifold cohomology agrees as a Hodge structure with the cohomology of any crepant resolution, should one exist.
As mentioned above, our strategy to construct $V^{(n)}$ hinges on finding a crepant resolution $\widetilde{X}$ of the hypersurface $X$.
This crepant resolution will be the $W^{(n - 1)}$ in \cref{intro/largebetti}.
So the $V$ and $W$ in \cref{intro/largeindex,intro/largebetti} respectively are related by
\[          
    V^{(n)} = \frac{W^{(n - 1)} \times E}{\mu_m}.
\]

We propose the following conjecture, which follows from the analogous conjecture of Esser, Totaro, and Wang due to the equality between orbifold Betti numbers and Betti numbers of a crepant resolution proven by Yasuda.
\begin{conjecture}
\label{intro/smoothbetti}
    $W^{(n)}$ has the largest sum of Betti numbers of any smooth, projective Calabi-Yau $n$-fold.
    When $n$ is odd, it also has the smallest negative Euler characteristic among such varieties.
\end{conjecture}

To prove our two main theorems, we are reduced to finding a $\mu_m$-equivariant, crepant, projective resolution of the weighted projective hypersurface $X$ we discussed in the construction of $V^{(n)}$.
As $X$ is \emph{quasi-smooth} in $\P$ (i.e. its affine cone is smooth away from the origin), and as the $\mu_m$ action is by toric automorphisms, we will show that to find an equivariant, crepant, projective resolution of $X$, it suffices to find a toric, crepant, projective resolution of $\P$.
Finding such a crepant resolution of $\P$ is the main technical work of this paper.

There is no guarantee that any crepant resolution of $\P$ exists.
Indeed, the singularities of a weighted projective space are of the form $\C^d/\mu_k$ where $\mu_k$ acts diagonally.
These do not always admit crepant resolutions.
For example, as discussed in \cite{nocrepant}, $\C^d/\mu_2$ is resolved by blowing up the origin, which leads to a single exceptional divisor with discrepancy $\frac{d}{2} - 1$.
This is therefore Gorenstein for $d$ even and terminal for $d > 2$.
Varieties which are $\Q$-factorial and terminal do not admit (non-trivial) crepant resolutions, so for instance $\C^{4}/\mu_2$ is a Gorenstein quotient singularity with no crepant resolution.

The strongest general result we are aware of is by Dais, Haase, and Ziegler in \cite{lci}, stating that toric local complete intersection (lci) singularities admit projective, crepant resolutions.
This is inapplicable in our setting, as the singularities we resolve are in general not lci.

The paper is structured as follows:
\begin{itemize}
    \item \cref{section/background} describes some of the relevant notation and background on Sylvester's sequence and lattice polytopes.
    \item \cref{section/toric} provides background information on toric varieties, and describes in combinatorial terms what is needed to provide a toric, crepant, projective resolution of a weighted projective space.
    \item \cref{section/triangles} describes a solution to this combinatorial problem for the weighted projective spaces in question, and thereby comes to toric, crepant, projective resolutions of them.
    \item \cref{section/hypersurface} uses these resolutions to resolve the hypersurfaces in question and provide examples of smooth, projective Calabi-Yau varieties of large index and Betti numbers, proving our main theorems.
\end{itemize}

\textbf{Acknowledgements.}
This work was supported by NSF grants DMS-2054553 and DMS-2136090.
I would like to thank Valery Alexeev, Louis Esser, Joaqu\'in Moraga, and Burt Totaro for the invaluable discussions.

    \section{Notation and background}
        \label{section/background}
        \subsection{Sylvester's sequence}
\label{section/background/sylvester}

\begin{definition}[A000058 on \cite{oeis}]
\hypertarget{background/sylvester}
    Sylvester's sequence $s_n$ is defined recursively via
    \begin{align*}
        s_0 &= 2\\
        s_{n + 1} &= s_0 s_1 \cdots s_n + 1.
    \end{align*}
\end{definition}

The recurrence relationship can be rewritten as
\[
    \sum_{i = 0}^{n} \frac{1}{s_i} = 1 - \frac{1}{s_{n + 1} - 1},
\]
so it follows that
\[
    \sum_{i = 0}^{\infty} \frac{1}{s_i} = 1,
\]
as $s_n \to \infty$.
This convergence is extremely fast, as Sylvester's sequence is doubly exponential.  In fact, $s_n > 2^{2^{n - 1}}$.
This is the fastest a sum of unit fractions can converge to $1$.
That is, any sum of $n$ unit fractions which is less than $1$ is at most $\sum_{i = 0}^{n - 1} \frac{1}{s_i}$, as proven in \cite{sylvester}.

These properties have been used in a number of recent results by Esser, Totaro, and Wang in finding conjectural extremal values for various numerical invariants of varieties in all dimensions, such as the index and sum of Betti numbers described in \cref{section/intro}.

\subsection{Lattice polytopes}
\label{section/background/polytope}

The combinatorial problem we must solve to find crepant resolutions of weighted projective spaces comes down to triangulating lattice polytopes (a \emph{lattice point} in $\R^n$ is a point in $\Z^n$).
We refer to \cite{unimodular} by Haase, Paffenholz, Piechnik, and Santos for more details on lattice polytopes.
By convention, all of our polytopes are convex.

\begin{definition}
    A \emph{(lattice) polytope} in $\R^n$ is the convex hull of a finite set of (lattice) points in $\R^n$.
    A \emph{(lattice) simplex} in $\R^n$ is the convex hull of a nonempty affinely independent set of (lattice) points in $\R^n$.
\end{definition}

A hyperplane $H$ in $\R^n$ splits $\R^n$ into two closed half-spaces.
If a polytope $P$ is wholly contained in one of these half-spaces, we say $H \cap P$ is a \emph{face} of $P$, which is itself a polytope.
A \emph{facet} of $P$ is a face of codimension $1$.
By definition, the empty set and $P$ are both faces of $P$.

A polytope can be defined either as a convex hull of finitely many vertices or as an intersection of finitely many closed half-spaces.
For example,
\[
    \conv(e_1, \dots, e_n, 0) = \curly{x \in \R^n \; \bigg| \; 1 - \sum_{i = 1}^{n} x_i \geq 0} \cap \bigcap_{i = 1}^n \curly{x \in \R^n \; \bigg| \; x_i \geq 0},
\]
where $e_i$ denotes the $i^{th}$ standard basis vector in $\R^n$.

Moreover, an $n$-dimensional simplex in $\R^n$ with the origin in its interior can be described as the convex hull of an affine basis $\{v_0, \dots, v_n\}$ of $\R^n$, or as the intersection
\[
    \bigcap_{i = 0}^{n} \{x \in \R^n \mid \Angle{u_i, x} + 1 \geq 0\},
\]
for an affine basis $\{u_0, \dots, u_n\}$ of $\R^n$ whose convex hull contains the origin.

We can swap between these representations using \emph{polar duality}.
Our reference for this fact is \cite[Theorem 2.11]{polytope}.
\begin{definition}
\hypertarget{background/polardual}
    Let $P$ be a polytope in $\R^n$.
    The \emph{polar dual} of $P$ is
    \[
        \widecheck{P} = \{y \in \R^n \mid \Angle{x, y} + 1 \geq 0 \text{ for all } x \in P\}.
    \]
\end{definition}

Note that in \cite{polytope}, Ziegler uses a slightly different definition of (and notation for) the polar dual of a polytope.
We choose the definition written here as it is the one used in SageMath \cite{sagemath}, which has been an invaluable computational tool.

\begin{proposition}[Polar duality]
\label{background/polarduality}
    Let $\sigma$ be an $n$-dimensional simplex in $\R^n$ containing $0$ in its interior.
    \begin{itemize}
        \item $\widecheck{\sigma}$ is an $n$-dimensional simplex in $\R^n$ containing $0$ in its interior.
        \item If
        \[
            \sigma = \conv(v_0, \dots, v_n) = \bigcap_{i = 0}^{n} \{x \in \R^n \mid \Angle{u_i, x} + 1 \geq 0\},
        \]
        then
        \[
            \widecheck{\sigma} = \conv(u_0, \dots, u_n) = \bigcap_{i = 0}^{n} \{x \in \R^n \mid \Angle{v_i, x} + 1 \geq 0\}.
        \]
    \end{itemize}
\end{proposition}

\begin{definition}
    Let $P$ be a (lattice) polytope.
    A \emph{(lattice) subdivision} $\mathcal S$ of $P$ is a finite set of (lattice) polytopes with the following properties:
    \begin{itemize}
        \item Every face of every member of $\mathcal S$ is a member of $\mathcal S$.
        \item The intersection of any two elements of $\mathcal S$ is a common face of both of them.
        \item The union of the polytopes in $\mathcal S$ is $P$.
    \end{itemize}
\end{definition}

A polytope $Q \in \mathcal S$ is called a \emph{cell} if $\dim(Q) = \dim(P)$.
If every nonempty element of $\mathcal S$ is a simplex, we say $\mathcal S$ is a \emph{triangulation}.

For a $k$-dimensional polytope $P \subseteq \R^n$ we define its \emph{normalized volume} to be
\hypertarget{background/nvol}
\[
    \nvol(P) = k! \, \vol(P),
\]
where $\vol(P)$ is the volume of $P$ in the $k$-dimensional subspace $\mathrm{span}(P)$.
This allows for convenient comparison of the volume of simplices in different dimensions.
We compute the normalized volume of a simplex $\sigma = \conv(v_0, \dots, v_n)$ in $\R^n$ as
\[
    \nvol(\sigma) = \pipe{\det{
            \begin{pmatrix}
                v_0 & \dots & v_n\\
                -1 & \dots & -1
            \end{pmatrix}
        }
    }
    ,
\]
where we view the $v_i$ as column vectors.

\begin{convention}
    \emph{Polytope} means lattice polytope and \emph{subdivision} means lattice subdivision.
\end{convention}

\begin{definition}
\hypertarget{background/unimodular}
    We say a simplex $\sigma$ is \emph{unimodular} if $\nvol(\sigma) = 1$.
    Furthermore, a triangulation $\mathcal T$ is called \emph{unimodular} if its cells are unimodular.
\end{definition}

If a simplex $\sigma$ admits a unimodular triangulation $\mathcal T$ (which is not guaranteed past dimension 2), $\nvol(\sigma)$ is the number of cells in $\mathcal T$.

\begin{definition}
\hypertarget{background/regular}
    A subdivision $\mathcal S$ of an $n$-dimensional polytope $P$ in $\R^n$ is called \emph{regular} if there is a continuous, piecewise linear function $\psi: P \longrightarrow \R$ so that $\psi(P \cap \Z^n) \subseteq \Q$ and is \emph{strictly convex} with respect to $\mathcal S$.
    That is, $\psi$ is convex and the domains of linearity of $\psi$ are precisely the cells of $\mathcal S$.
    (The graph of $\psi$ can be viewed as a folding of the subdivision $\mathcal S$ into a strictly convex shape.)
\end{definition}
Note that some author's definition of regularity imposes $\psi(P \cap \Z^{n}) \subseteq \Z$ rather than $\Q$.
This is equivalent to the definition we give.
Indeed, the polytopes we are working with are compact and hence $P \cap \Z^n$ is finite, so such a $\psi$ with $\psi(P \cap \Z^{n}) \subseteq \Q$ actually has image in $\frac{1}{N} \Z$ for a positive integer $N$.
We may then multiply $\psi$ by $N$ to land in $\Z$.

\begin{definition}
\hypertarget{background/starshaped}
    A subdivision $\mathcal S$ is called \emph{star-shaped} (at the origin) if every cell contains the origin as a vertex.
\end{definition}

	\section{Background on resolving toric varieties}
        \label{section/toric}
        We describe in this section our strategy to find toric, crepant, projective resolutions of weighted projective spaces.
Our main reference for facts about toric varieties is the short text \cite{toric}, whose notation we borrow.
We recall some of the basic facts here.

A toric variety is defined by a lattice $N$ and a strongly convex rational polyhedral fan $\Delta$ in $N \otimes \R$.
We denote this variety by $X(\Delta)$, which comes equipped with a natural action and equivariant open immersion of the algebraic torus $T_{X(\Delta)} = T_N = N \otimes \C^\times$.

Key geometric properties of $X(\Delta)$ correspond to combinatorial properties of the fan $\Delta$, which encodes the stratification of $X(\Delta)$ into orbit closures.
For example, we will crucially use the following characterization of smoothness.
\begin{proposition}\cite[\S 2.1]{toric}
\label{toric/smoothness}
    A toric variety $X(\Delta)$ is smooth if and only if every maximal cone $\sigma \in \Delta$ is generated by part of an integral basis of the lattice $N$.
\end{proposition}

Similarly, torus invariant divisors can be fully encoded combinatorially.
Let $\Delta(1)$ be the set of one dimensional cones, i.e. rays, in $\Delta$.
Then for every $\rho \in \Delta(1)$, there is a corresponding orbit closure $D_{\rho} = V(\rho)$ which is a torus invariant, irreducible Weil divisor.
Every torus invariant Weil divisor is of the form $\sum_{\rho \in \Delta(1)} a_\rho D_\rho$ for $a_{\rho} \in \Z$.
For example, the canonical divisor of any toric variety $X(\Delta)$ is given by
\[
    K_{X(\Delta)} = -\sum_{\rho \in \Delta(1)} D_{\rho}.
\]

Torus invariant ($\Q$-)Cartier divisors correspond to continuous functions $\psi: |\Delta| \longrightarrow \R$ which restrict to integral (rational) linear functions on every cone in $\Delta$.
Here, $|\Delta| = \bigcup_{\sigma \in \Delta} \sigma$ is the support of the fan $\Delta$.
If a torus invariant Weil divisor $D = \sum_{\rho \in \Delta(1)} a_{\rho} D_{\rho}$ is Cartier, the associated piecewise linear continuous function $\psi_D: |\Delta| \longrightarrow \R$ is called the support function of $D$, and is uniquely determined by the condition $\psi_D(n_{\rho}) = -a_{\rho}$, where $n_{\rho}$ is the minimal nonzero lattice point on the ray $\rho$.

Important properties of a torus invariant Cartier divisor can be determined from its support function.
For our purposes, we will primarily need the following characterization of ampleness.
\begin{proposition}\cite[\S 3.4]{toric}
\label{toric/divisors}
    Let $X(\Delta)$ be a toric variety so that $|\Delta| = N \otimes \R$.
    A torus invariant Cartier divisor $D$ on $X(\Delta)$ is ample if and only if its support function $\psi$ is strictly convex, i.e. it is convex and for distinct maximal cones $\sigma, \sigma' \in \Delta$, the linear functions extending $\psi|_{\sigma}$ and $\psi|_{\sigma'}$ are distinct.
\end{proposition}

The construction of toric varieties is functorial.
Let $(N, \Delta)$ and $(N', \Delta')$ be lattices equipped with fans.
Given a $\Z$-linear morphism $f: N \longrightarrow N'$ so that for every $\sigma \in \Delta$ there is a $\sigma' \in \Delta'$ satisfying $f(\sigma) \subseteq \sigma'$, there is a morphism $\phi: X(\Delta) \longrightarrow X(\Delta')$ which is equivariant with respect to the torus actions.
On tori, $\phi$ is given by $f \otimes \id_{\C^\times} : N \otimes \C^\times \longrightarrow N' \otimes \C^\times$.
             
\begin{proposition}
\label{toric/birational}
    $f: N \longrightarrow N'$ is an isomorphism if and only if $\phi: X(\Delta) \longrightarrow X(\Delta')$ is birational.
\end{proposition}
\begin{proof}
    If $f$ is an isomorphism then $\phi|_{T_N} = f \otimes \id : N \otimes \C^{\times} \longrightarrow N' \otimes \C^{\times}$ is an isomorphism, so $\phi$ is birational.

    Conversely, suppose $\phi$ is birational.
    Then $\phi|_{T_N}: T_N \longrightarrow T_{N'}$ is a birational homomorphism of algebraic tori.
    As such, it is an isomorphism.
    As $N \mapsto N \otimes \C^{\times}$ induces an equivalence of categories between lattices and algebraic tori, it follows that $f$ is an isomorphism.
\end{proof}

\begin{proposition}
\label{toric/pullbacksupport}
    Given $\phi: X(\Delta) \longrightarrow X(\Delta')$ arising from $f: N \longrightarrow N'$ and a torus invariant Cartier divisor $D'$ on $X(\Delta')$ with support function $\psi_{D'}$, we have $\psi_{\phi^* D'} = \psi_{D'} \circ f$.
\end{proposition}
\begin{proof}
    Let $M'$ be the dual of $N'$.
    For $\sigma' \in \Delta'$, let $-u \in M'$ be the linear function whose restriction to $\sigma'$ agrees with $\psi_{D'}$, so that $D'|_{U_{\sigma'}} = \div(\chi^{u})$.
    Let $\sigma \in \Delta$ be so that $f(\sigma) \subseteq \sigma'$.
    Then $\phi$ restricts to a morphism of affine toric varieties $U_{\sigma} \longrightarrow U_{\sigma'}$ defined on coordinate rings $A_{\sigma'} \longrightarrow A_{\sigma}$ by taking a monomial $\chi^{m}$ to $\chi^{m \circ f}$.
    Hence, $\phi^*(D')|_{U_{\sigma}} = \div(\chi^{u \circ f})$.
    As such, we have shown that $\psi_{\phi^* D'}$ agrees with $\psi_{D'} \circ f$ on any cone $\sigma' \in \Delta'$.
    Hence, they agree on the whole support of $\Delta'$.
\end{proof}

A common way to construct toric varieties is to take a polytope $P$ with the origin in its relative interior and to form the fan $\Delta_P$ consisting of cones generated by all the proper faces of $P$.
For example, in the lattice $\Z^{n + 1}/(1, \dots, 1)$, the polytope $\conv(\overline{e_0}, \dots, \overline{e_n})$ (where $\overline{e_i}$ is the image of the $i^{th}$ standard basis vector in the quotient $\Z^{n + 1}/(1, \dots, 1)$) determines a very familiar toric variety, namely $\P^n$.
More generally, given positive integer weights $a_0, \dots, a_n$, the polytope $\conv(\overline{e_0}, \dots, \overline{e_n})$ in the lattice $\Z^{n + 1}/(a_0, \dots, a_n)$ defines the weighted projective space $\P(a_0, \dots, a_n)$.
(If $a_n = 1$, we can project to the first $n$ coordinates and work in the lattice $\Z^n$ with the polytope $P = \conv(e_0, \dots, e_{n - 1}, -\sum a_i e_i)$.)
More generally still, if we have a \hyperlink{background/starshaped}{star-shaped} subdivision $\mathcal S$ of a polytope $P$, we can form the fan $\Delta_{\mathcal S}$ consisting of cones generated by all the elements of $\mathcal S$.
For the remainder of this section, $\mathcal S$ is a star-shaped subdivision of $P$.

The identity morphism $\id: N \longrightarrow N$ sends cones in $\Delta_{\mathcal S}$ to $\Delta_{P}$ and hence defines a morphism of toric varieties $\pi: X(\Delta_{\mathcal S}) \longrightarrow X(\Delta_{P})$.
It is automatic that $\pi$ is equivariant, and by \cref{toric/birational} it's also birational.
We will show below that $\pi$ is crepant under some hypotheses, and determine conditions for when $\pi$ is a projective resolution of singularities.
As our goal is to resolve weighted projective spaces, we will focus on the case of interest where $P$ is the simplex $\conv(\overline{e_0}, \dots, \overline{e_n})$ in $\Z^{n + 1}/(a_0, \dots, a_n)$.
Let $\P = \P(a_0, \dots, a_n) = X(\Delta_P)$.

First, we show that $\pi: X(\Delta_\mathcal S) \longrightarrow \P$ is crepant.
Let $\psi$ be the support function of $K_\P$, which is uniquely determined by $\psi(\overline{e_i}) = 1$ for all $i$.
Let's suppose $\psi$ is piecewise $\Z$-linear, i.e. that $K_{\P}$ is Cartier.
This is not a substantial assumption to make, as any normal variety with a $\Q$-Cartier canonical divisor admitting a crepant resolution must actually have a Cartier canonical divisor, and this assumption is satisfied in the case of the weighted projective spaces we study in \cref{section/triangles}.
As $\pi$ is defined by $\id: N \longrightarrow N$, we have by \cref{toric/pullbacksupport} that $\pi^*(K_{\P})$ has support function $\psi$.
The support function of $K_{X(\Delta_\mathcal S)}$ is defined by $\psi_{K_{X(\Delta_\mathcal S)}}(n_\rho) = 1$ for every ray $\rho$ of $\Delta_\mathcal S$, where $n_\rho$ is the minimal lattice point on $\rho$.
Any such ray $\rho$ is generated by some lattice point $m \in \partial P$.
As such, $\psi(m) \leq 1$.
Because $\psi$ is piecewise $\Z$-linear, $\psi(n_{\rho}) \in \Z$.
Furthermore, $0 < \psi(n_{\rho}) \leq \psi(m) \leq 1$.
Hence, $\psi(n_{\rho}) = 1$ for every ray $\rho$ of $\Delta_{\mathcal S}$.
As such, $\psi$ equals the support function of the canonical divisor $\psi_{K_{X(\Delta_{\mathcal S})}}$ which proves that $\pi^*(K_{X(\Delta_P)}) = K_{X(\Delta_{\mathcal S})}$.
That is, $\pi$ is crepant.

$X(\Delta_{\mathcal S})$ is projective if and only if it has an ample $\Q$-Cartier divisor.
By \cref{toric/divisors} above, it suffices to find a piecewise $\Q$-linear continuous function $|\Delta_{\mathcal S}| \longrightarrow \R$ which is strictly convex with respect to $\mathcal S$.
That is, it suffices to show that $\mathcal S$ is a \hyperlink{background/regular}{regular subdivision} of $P$ as defined in \cref{section/background/polytope}.

Finally, 
$X(\Delta_{\mathcal S})$ is smooth if and only if the nonzero vertices of every cell $\sigma \in \mathcal S$ form an integral basis of the lattice, i.e. are \hyperlink{background/unimodular}{unimodular simplices}, as defined in \cref{section/background/polytope}.

Altogether, we have determined that to find a toric, crepant, projective resolution of a weighted projective space $\P$, we must find a regular, unimodular, star-shaped triangulation $\mathcal T$ of the simplex $P$ corresponding to $\P$ described above.

	\section{Triangles}
        \label{section/triangles}
        Let us proceed to find such a triangulation, as described in \cref{section/toric}.
We define
\begin{align*}
    d_1^{(n)} &= 2 s_{n - 1} - 2\\
    d_2^{(n)} &= s_{n} - 1
\end{align*}
and 
\begin{align*}
    w_1^{(n)} &= \parens{-\frac{d_1^{(n)}}{s_0}, \dots, -\frac{d_1^{(n)}}{s_{n - 2}}, -1}\\
    w_2^{(n)} &= \parens{-\frac{d_2^{(n)}}{s_0}, \dots, -\frac{d_2^{(n)}}{s_{n - 1}}}.
\end{align*}
With this, we define
\begin{align*}
    P_1^{(n)} &= \conv\parens{e_0^{{(n)}}, \dots, e_{n-1}^{(n)}, w_1^{(n)}}\\
    P_2^{(n)} &= \conv\parens{e_0^{{(n)}}, \dots, e_{n-1}^{(n)}, w_2^{(n)}},
\end{align*}
where $e_i^{{(n)}}$ is the $i^{th}$ standard basis vector in $\R^n$ (indexed beginning with $0$) and $\conv$ denotes the convex hull.
When it is clear from context, the superscript $(n)$ will be omitted for brevity.

These are both simplices in $\R^n$, and correspond as in \cref{section/toric} to the weighted projective spaces
\begin{align*}
    \P_1^{(n)} &= \P\parens{d_1^{(n)} / s_0, \dots, d_1^{(n)}/s_{n - 2}, 1, 1}\\
    \P_2^{(n)} &= \P\parens{d_2^{(n)} / s_0, \dots, d_2^{(n)} / s_{n - 1}, 1}.
\end{align*}

To use the strategy from \cref{section/toric} to crepantly resolve these weighted projective spaces, we must first verify that the canonical divisors $K_{\P_i^{(n)}}$ are Cartier for $i = 1, 2$ and any $n \geq 1$.
By \cite[Corollary 4A.5]{introwps}, a torus invariant divisor $D = \sum c_i D_i$ in $\P(a_0, \dots, a_n)$ is Cartier if $\mathrm{lcm}(a_0, \dots, a_n) | \sum c_i a_i$, where $D_i = \{x_i = 0\}$ is the torus invariant divisor corresponding to the ray generated by $\overline{e_i}$.
The converse is true supposing $\P(a_0, \dots, a_n)$ is \hyperlink{hypersurface/wellformed}{\emph{well-formed}}, as will be defined in \cref{section/hypersurface}.
For $D = K_{\P}$, this means that $K_{\P}$ is Cartier if $a_i \mid \sum_j a_j$ for all $i$.
This holds for both weighted projective spaces in question.
We will use the relation:
\[
    \sum_{j = 0}^{k} \frac{1}{s_j} = 1 - \frac{1}{s_{k + 1} - 1}.
\]

For $i = 1$, we compute the sum of the weights to be
\begin{align*}
    \sum_{j = 0}^{n - 2} \frac{d_1^{(n)}}{s_j} + 1 + 1 &= (2 s_{n - 1} - 2)\parens{1 - \frac{1}{s_{n - 1} - 1}} + 2\\
    &= 2 s_{n - 1} - 2\\
    &= d_1^{(n)}.
\end{align*}
All of the weights divide $d_1^{(n)}$, so we deduce that $K_{\P_1^{(n)}}$ is Cartier for all $n$.

For $i = 2$, the sum of the weights is
\begin{align*}
    \sum_{j = 0}^{n - 1} \frac{d_2^{(n)}}{s_j} + 1 &= (s_{n} - 1)\parens{1 - \frac{1}{s_n - 1}} + 1\\
    &= s_n - 1\\
    &= d_2^{(n)}.
\end{align*}
All of the weights divide $d_2^{(n)}$, so $K_{\P_2^{(n)}}$ is Cartier for all $n$.

The sources of the large index examples in \cite{largeindex} are hypersurfaces in $\P_1$.
We seek to find a toric, crepant, projective resolution of $\P_1$.
To do so, as discussed in \cref{section/toric}, is equivalent to finding a regular, unimodular, star-shaped triangulation of $P_1$.
\begin{proposition}
\label{triangles/triangulating1}
    For all $n$, there is a regular, unimodular, star-shaped triangulation of $P_1^{(n)}$.
\end{proposition}
In the following section, we reduce this to finding a regular, unimodular, star-shaped triangulation of $P_2$.

\subsection{Reducing to $P_2$}
\label{section/triangles/reduction}

Notice that $w_1^{{(n+1)}} = (2w_2^{{(n)}}, -1)$.
From this, we can show that the family $P_2^{(n)}$ sits inside the family $P_1^{(n + 1)}$ as the hyperplane section setting $x_{n} = 0$.
See \cref{triangles/figure/2inside1}.
\begin{figure}
    \begin{minipage}{\textwidth}
        \centering
        \begin{tikzpicture}[scale = 2]
    \fill[color=pink] (1, 0) -- (-2, -1) -- (0,1) -- cycle {};

    \draw[color=black, thick] (1, 0) -- (-1, 0);
    \node[below] at (0, 0) {$P_2^{(1)} \times \{0\}$};
    \draw[color=black, thick, opacity=0.6] (1, 0) -- (2, 0);
    \node[black, above] at (1.5, 0) {$\{x_1 = 0\}$};
    \draw[color=black, thick, opacity=0.6] (-1, 0) -- (-3, 0);

    \node[above] at (0, 1) {$e_1^{(2)}$};
    \node[below] at (-2, -1){$w_1^{(2)}$};
\end{tikzpicture}
        \caption*{(a) $P_1^{(2)}$ with its cross section $P_2^{(1)} \times \{0\} = P_1^{(2)} \cap \{x_1 = 0\}$.}
    \end{minipage}

    \vspace{1cm}

    \begin{minipage}{\textwidth}
        \centering
        \begin{tikzpicture}%
    [x={(-0.094095cm, -0.098890cm)},
    y={(0.995563cm, -0.009478cm)},
    z={(0.000131cm, 0.995053cm)},
    scale=2,
    back/.style={very thick, black},
    edge/.style={color=red!95!red, thick, opacity=0.3},
    facet/.style={fill=blue!95!black,fill opacity=1},
    vertex/.style={inner sep=1pt,circle,draw=black!25!black,fill=black!75!black,thick}]
    %
    %
    \definecolor{rightpink}{HTML}{ff5e5e}
    \definecolor{leftpink}{HTML}{ff7575}

    \fill[color=black, opacity=0.2] (-5, -3.5, 0)  -- (3, 2.5, 0) -- (-4, 2.5, 0);

    \fill[color=leftpink, opacity=0.7] (1.00000, 1.00000, 0.00000) -- (-6.00000, -4.00000, -1.00000) -- (0.00000, 0.00000, 1.00000) -- cycle {};
    \fill[color=rightpink, opacity=0.7] (1.00000, 1.00000, 0.00000) -- (0.00000, 0.00000, 1.00000) -- (-3.00000, 1.00000, 0.00000) -- cycle {};


    \fill[color=black, opacity=0.2] (-5, -3.5, 0) -- (0, -4.25, 0) -- (3, 2.5, 0);

    \draw
        (-3, -2, 0) edge[line width=1, opacity=0.2, dashed] (-3, 1, 0)
        (-3, -2, 0) edge[line width=1, opacity=0.2] (1, 1, 0)
        (-3, 1, 0) edge[line width =1, opacity=0.2] (1, 1, 0);

    \node[color=black, above] at (0, 0, 1) {$e_2^{(3)}$};
    \node[color=black, below] at (-6, -4, -1) {$w_1^{(3)}$};

    \node[color=black, above] at (-4, 2, 0) {$\{x_2 = 0\}$};

    \node[color=black, above] at (0, 0, 0.28) {$P_2^{(2)} \times \{0\}$};
\end{tikzpicture}
        \caption*{(b) $P_1^{(3)}$ with its cross section $P_2^{(2)} \times \{0\} = P_1^{(3)} \cap \{x_2 = 0\}$.}
    \end{minipage}
    \caption{}
    \label{triangles/figure/2inside1}
\end{figure}
More precisely, we can show the following.
\begin{lemma}
\label{triangles/2inside1}
    $P_1^{{(n + 1)}} = \conv\parens{P_2^{{(n)}} \times \{0\}, e_n^{{(n + 1)}}, w_1^{{(n + 1)}}}$.
\end{lemma}
\begin{proof}
    We first show $\conv(P_2^{{(n)}} \times \{0\}, e_n^{{(n + 1)}}, w_1^{{(n + 1)}}) \subseteq P_1^{(n + 1)}$.
    Indeed, $e_n^{(n + 1)}$ and $w_1^{(n + 1)}$ are in $P_1^{(n + 1)}$ by definition.
    Recall that $P_2^{(n)} \times \{0\}$ is the convex hull of the first $n$ standard basis vectors in $\R^{n + 1}$ and $(w_2^{(n)}, 0)$.
    The standard basis vectors are all, by definition, in $P_1^{(n + 1)}$.
    As for $(w_2^{(n)}, 0)$, we may realize it as a convex combination of the vertices of $P_1^{(n + 1)}$ via
    \[
	\parens{w_2^{(n)}, 0} = \frac{1}{2} w_1^{(n + 1)} + \frac{1}{2} e_n^{(n + 1)},
    \]
    showing that it is also in $P_1^{(n + 1)}$.
    As $P_1^{(n + 1)}$ is convex, the containment $$\conv\parens{P_2^{{(n)}} \times \{0\}, e_n^{{(n + 1)}}, w_1^{{(n + 1)}}} \subseteq P_1^{(n + 1)}$$ is proven.

    As we have this containment, to show that these polytopes are equal to each other it suffices to show that they have the same volume.
    Indeed, they have the same dimension so it suffices to compare the \hyperlink{background/nvol}{normalized volume}.
    The normalized volume of $P_1^{(n + 1)}$ is
    \[
	\nvol\parens{P_1^{(n + 1)}}
        =
        \pipe{\det
            \begin{pmatrix}
                & & & \vline & \vline\\
                & & & \vline & \vline\\
                & I_{n + 1} & & \vline & w_1^{(n + 1)}\\
                & & & \vline & \vline\\
                & & & \vline & \vline\\
                \hline
                -1 & \dots & -1 & & -1
            \end{pmatrix}
        }.
    \]
    Recall that $w_1^{(n + 1)} = (2 w_2^{(n)}, -1)$.
    We may then rewrite the above as
    \[
        \pipe{\det
            \begin{pmatrix}
                & & & \vline & \vline\\
                & & & \vline & \vline\\
                & I_{n + 1} & & \vline & 2w_2^{(n)}\\
                & & & \vline & \vline\\
                & & & \vline & -1\\
                \hline
                -1 & \dots & -1 & & -1
            \end{pmatrix}
        }
        ,
    \]
    which we compute as
    \begin{align*}
        \pipe{\det
            \begin{pmatrix}
                & & & \vline & \vline\\
                & & & \vline & \vline\\
                & I_{n + 1} & & \vline & 2w_2^{(n)}\\
                & & & \vline & \vline\\
                & & & \vline & -1\\
                \hline
                -1 & \dots & -1 & & -1
            \end{pmatrix}
        }
        &=
        \pipe{\det
            \begin{pmatrix}
                & & & \vline & \vline\\
                & & & \vline & \vline\\
                & I_{n + 1} & & \vline & 2w_2^{(n)}\\
                & & & \vline & \vline\\
                & & & \vline & 0\\
                \hline
                -1 & \dots & -1 & & -2
            \end{pmatrix}
        }
        \\
        &=
        2 \pipe{\det
            \begin{pmatrix}
                & & & \vline & \vline\\
                & & & \vline & \vline\\
                & I_{n + 1} & & \vline & w_2^{(n)}\\
                & & & \vline & \vline\\
                & & & \vline & 0\\
                \hline
                -1 & \dots & -1 & & -1
            \end{pmatrix}
        }
        \\
        &=
	    2 \parens{\nvol\parens{P_2^{(n)}}}.
    \end{align*}

    Observe that $\conv\parens{P_2^{(n)} \times \{0\}, e_n^{(n + 1)}, w_1^{(n + 1)}}$ contains
    \[
	\conv\parens{P_2^{(n)} \times \{0\}, e_n^{(n + 1)}} \cup \conv\parens{P_2^{(n)} \times \{0\}, w_1^{(n + 1)}}.
    \]
    To compute the normalized volume of this union, we need only sum the normalized volumes of $\conv\parens{P_2^{{(n)}} \times \{0\}, e_n^{(n + 1)}}$ and $\conv\parens{P_2^{{(n)}} \times \{0\}, w_1^{(n + 1)}}$, as their intersection is $P_2^{(n)} \times \{0\}$ which has $0$ volume in $\R^{n + 1}$.
    To compute these two volumes, we appeal to the following lemma.

    \begin{lemma}
    \label{triangles/conevolume}
        Let $P = \conv(v_0, \dots, v_n)$ be a simplex in $\R^n$.
        Let $w \in \R^{n + 1}$ have last component $\pm 1$.
	Then $\nvol\parens{\conv(P \times \{0\}, w)} = \nvol\parens{P}$.
    \end{lemma}
    \begin{proof}
        Write $w = (w', \pm1)$.
	Then $\nvol\parens{\conv\parens{P \times \{0\}, w}}$ is given by the determinant
        \[
            \pipe{\det
                \begin{pmatrix}
                    v_0 & \dots & v_n & w'\\
                    0 & \dots & 0 & \pm 1\\
                    -1 & \dots & -1 & -1
                \end{pmatrix}
            }
            .
        \]
        By minor expansion, this determinant equals
        \[
            \pipe{\det
                \begin{pmatrix}
                    v_0 & \dots & v_n\\
                    -1 & \dots & -1
                \end{pmatrix}
            }
            ,
        \]
        which is precisely $\nvol(P)$.
    \end{proof}

    Using this, we compute
    \begin{align*}
	\nvol\parens{P_2^{(n)}} &= \nvol\parens{\conv\parens{P_2^{(n)} \times \{0\}, e_n^{(n + 1)}}}\\
	&= \nvol\parens{\conv\parens{P_2^{(n)} \times \{0\}, w_1^{(n + 1)}}}.
    \end{align*}
    We therefore know 
    \begin{align*}
        \nvol\parens{\conv\parens{P_2^{{(n)}} \times \{0\}, e_n^{{(n + 1)}}, w_1^{{(n + 1)}}}} &\geq \nvol\parens{\conv\parens{P_2^{{(n)}} \times \{0\}, e_n^{{(n + 1)}}}}\\
        &+ \nvol\parens{\conv\parens{P_2^{{(n)}} \times \{0\}, w_1^{{(n + 1)}}}}\\
	&= 2 \parens{\nvol\parens{P_2^{(n)}}}\\
        &= \nvol\parens{P_1^{(n + 1)}}\\
        &\geq \nvol\parens{\conv\parens{P_2^{{(n)}} \times \{0\}, e_n^{{(n + 1)}}, w_1^{{(n + 1)}}}}.
    \end{align*}
    So we have shown
    \[
        P_1^{{(n + 1)}} = \conv\parens{P_2^{{(n)}} \times \{0\}, e_n^{{(n + 1)}}, w_1^{{(n + 1)}}}.
    \]
\end{proof}

From this we can extend subdivisions of $P_2^{(n)}$ to $P_1^{(n + 1)}$ by taking cones.
\begin{definition}
    Let $H$ be a hyperplane in $\R^{n}$ and let $P$ be a polytope in $H$ with a subdivision $\mathcal S$.
    For $z \in \Z^n - H$ we define a subdivision $\cone(z, \mathcal S)$ of $\conv(z, P)$ via $$
    \cone(z, \mathcal S) = \{\conv(z, \sigma) \mid \sigma \in \mathcal S\}.$$
\end{definition}
By \cref{triangles/conevolume}, we have the following.
\begin{lemma}
\label{triangles/coneunimodular}
    Let $P$ be a simplex in $\R^{n} \times \{0\} \subseteq \R^{n + 1}$ with a unimodular triangulation $\mathcal T$.
    Let $z \in \R^{n} \times \{\pm 1\}$.
    Then $\cone(z, \mathcal T)$ is a unimodular triangulation of $\conv(z, P)$.
\end{lemma}
Given a triangulation $\mathcal T$ of $P_2^{(n)}$, we form triangulations $\cone\parens{e_n^{(n + 1)}, \mathcal T}$ and $\cone\parens{w_1^{(n + 1)}, \mathcal T}$ of $\conv\parens{e_n^{(n + 1)}, P_2^{(n)} \times \{0\}}$ and $\conv\parens{w_1^{(n + 1)}, P_2^{(n)} \times \{0\}}$.
By \cref{triangles/coneunimodular}, these are both unimodular.
Using \cref{triangles/2inside1}, we may glue these along their intersection $P_2^{(n)} \times \{0\}$ to get a unimodular triangulation $\mathcal T'$ of $P_1^{(n + 1)}$.

We now seek to show that if $\mathcal T$ is regular, then $\mathcal T'$ is regular as well.
To do so, we use the following two lemmas.
\begin{lemma}
\label{triangles/hyperplaneregular}
    Let $H$ be a hyperplane of $\R^{n}$ and $P \subseteq H$ a polytope with a regular subdivision $\mathcal S$.
    Let $z$ be a lattice point not in $H$.
    Then $\cone(z, \mathcal S)$ is a regular subdivision of $\conv(z, P)$.
    See \cref{triangles/figure/hyperplaneregular}.
\end{lemma}
    \begin{figure}[thp]
        \begin{minipage}{\textwidth}
            \centering
            \begin{center}
    \begin{tikzpicture}[scale = 2]
        \draw[color=red, line width=3] (0, 0) -- (1, 0);
        \draw[color=blue, line width=3] (1, 0) -- (2, 0);
        \draw[color=magenta, line width=3] (2, 0) -- (3, 0);

        \filldraw[black] (0, 0) circle (1.5pt);
        \filldraw[black] (1, 0) circle (1.5pt);
        \filldraw[black] (2, 0) circle (1.5pt);
        \filldraw[black] (3, 0) circle (1.5pt);
    \end{tikzpicture}
\end{center}
            \caption*{(a) A polytope $P \subseteq H$ with a regular subdivision $\mathcal S$.}
        \end{minipage}
        \begin{minipage}{\textwidth}
            \centering
            \begin{center}
    \begin{tikzpicture}[scale = 2]
        \fill[color=red, opacity=0.4] (0, 0) -- (1, 0) -- (0,1) -- cycle {};
        \fill[color=blue, opacity=0.4] (1, 0) -- (2, 0) -- (0,1) -- cycle {};
        \fill[color=magenta, opacity=0.4] (2, 0) -- (3, 0) -- (0,1) -- cycle {};

        \draw[color=red, line width=3] (0, 0) -- (1, 0);
        \draw[color=blue, line width=3] (1, 0) -- (2, 0);
        \draw[color=magenta, line width=3] (2, 0) -- (3, 0);

        \draw[color=black, line width=1.5] (0, 0) -- (0, 1);
        \draw[color=black, line width=1.5] (1, 0) -- (0, 1);
        \draw[color=black, line width=1.5] (2, 0) -- (0, 1);
        \draw[color=black, line width=1.5] (3, 0) -- (0, 1);

        \filldraw[black] (0, 0) circle (1.5pt);
        \filldraw[black] (1, 0) circle (1.5pt);
        \filldraw[black] (2, 0) circle (1.5pt);
        \filldraw[black] (3, 0) circle (1.5pt);
        \filldraw[black] (0, 1) circle (1.5pt);

        \draw[color=black, line width=1.5, opacity=0.6] (-1, 0) -- (0, 0);
        \draw[color=black, line width=1.5, opacity=0.6] (3, 0) -- (4, 0);

        \node[above] at (4, 0) {$H$};
        \node[above] at (0, 1.05) {$z$};
    \end{tikzpicture}
\end{center}
            \caption*{(b) $\conv(z, P)$ with the subdivision $\cone(z, \mathcal S)$.}
        \end{minipage}
        \caption{}
        \label{triangles/figure/hyperplaneregular}
    \end{figure}
\begin{proof}
    Let $\psi: P \longrightarrow \R$ be a continuous, piecewise $\Q$-linear function which is strictly convex with respect to $\mathcal S$.
    That is, $\psi$ witnesses the regularity of $\mathcal S$.
    We extend $\psi$ to $\Psi: \conv(z, P) \longrightarrow \R$ by defining $\Psi|_P = \psi$, setting $\Psi(z)$ to be any rational number $\omega$, and extending $\psi$ linearly to $z$.
    That is, for any $y \in \conv(z, P)$ (other than $z$ itself) we consider the line connecting $z$ to $y$ and let $x$ be its intersection with $P$.
    Then $y = tx + (1 - t)z$ for a unique $t \in [0, 1]$, and we define $\Psi(y) = t \psi(x) + (1 - t) \omega$.
    See \cref{triangles/figure/hyperplaneregulargraph} for a visual of $\Psi$.

    \begin{figure}[thp]
        \begin{minipage}{0.5\textwidth}
            \centering
            \begin{center}
    \begin{tikzpicture}[scale = 2]
        \draw[color=red, line width=3] (0, 1) -- (1, 0);
        \draw[color=blue, line width=3] (1, 0) -- (2, 0);
        \draw[color=magenta, line width=3] (2, 0) -- (3, 1);

        \draw[color=red, line width=3] (0, -1) -- (1, -1);
        \draw[color=blue, line width=3] (1, -1) -- (2, -1);
        \draw[color=magenta, line width=3] (2, -1) -- (3, -1);

        \filldraw[black] (0, 1) circle (1.5pt);
        \filldraw[black] (1, 0) circle (1.5pt);
        \filldraw[black] (2, 0) circle (1.5pt);
        \filldraw[black] (3, 1) circle (1.5pt);

        \filldraw[black] (0, -1) circle (1.5pt);
        \filldraw[black] (1, -1) circle (1.5pt);
        \filldraw[black] (2, -1) circle (1.5pt);
        \filldraw[black] (3, -1) circle (1.5pt);

        \draw
            (0, 1) edge[line width=1, dashed] (0, -1)
            (1, 0) edge[line width=1, dashed] (1, -1)
            (2, 0) edge[line width=1, dashed] (2, -1)
            (3, 1) edge[line width=1, dashed] (3, -1);

        \node[below] at (1.5, -1.1) {$P$};

    \end{tikzpicture}
\end{center}
            \caption*{(a) The graph of $\psi$ lying over $P$.}
        \end{minipage}
        \begin{minipage}{0.5\textwidth}
            \centering
            \begin{center}
    \begin{tikzpicture}%
        [x={(0.760444cm, -0.251632cm)},
        y={(0.649404cm, 0.294756cm)},
    	z={(-0.000070cm, 0.921846cm)},
	    scale=2.500000,
	    back/.style={loosely dotted, thin},
	    edge/.style={color=blue!95!black, thick},
	    facet/.style={fill=blue!95!black,fill opacity=0.800000},
	    vertex/.style={inner sep=1pt,circle,draw=green!25!black,fill=green!75!black,thick}]
        %
        %

        \definecolor{topleft}{HTML}{880000}
        \definecolor{topcenter}{HTML}{1010fa}
        \definecolor{topright}{HTML}{f106f1}
        \definecolor{bottomleft}{HTML}{d00000}
        \definecolor{bottomcenter}{HTML}{0000d0}
        \definecolor{bottomright}{HTML}{d000d0}

        \fill[color=bottomleft, opacity=0.6] (1.00000, 0.00000, 0) -- (0.00000, 0.00000, 0) -- (0.00000, 1.00000, 0) -- cycle {};
        \fill[color=bottomcenter, opacity=0.6] (1.00000, 0.00000, 0) -- (2.00000, 0.00000, 0) -- (0.00000, 1.00000, 0) -- cycle {};
        \fill[color=bottomright, opacity=0.6] (2.00000, 0.00000, 0) -- (3.00000, 0.00000, 0) -- (0.00000, 1.00000, 0) -- cycle {};

        \draw
            (0, 0, 1.8) edge[line width=1, dashed] (0, 0, 0)
            (1, 0, 1) edge[line width=1, dashed] (1, 0, 0)
            (2, 0, 1) edge[line width=1, dashed] (2, 0, 0)
            (3, 0, 1.8) edge[line width=1, dashed] (3, 0, 0)
            (0, 1, 0.55) edge[line width=1, dashed] (0, 1, 0);
        
        \fill[color=topleft, opacity=0.6] (1.00000, 0.00000, 1.00000) -- (0.00000, 0.00000, 1.80000) -- (0.00000, 1.00000, 1.00000) -- cycle {};
        \fill[color=topcenter, opacity=0.6] (1.00000, 0.00000, 1.00000) -- (2.0000, 0.00000, 1.00000) -- (0.00000, 1.00000, 1.00000) -- cycle {};
        \fill[color=topright, opacity = 0.6] (2.00000, 0.00000, 1.00000) -- (3.00000, 0.00000, 1.80000) -- (0.00000, 1.00000, 1.00000) -- cycle {};

        \draw
            (0, 0, 0) edge[line width=2] (0, 1, 0)
            (1, 0, 0) edge[line width=2] (0, 1, 0)
            (2, 0, 0) edge[line width=2] (0, 1, 0)
            (3, 0, 0) edge[line width=2] (0, 1, 0)
            (0, 0, 0) edge[line width=4, color=red] (1, 0, 0)
            (1, 0, 0) edge[line width=4, color=blue] (2, 0, 0)
            (2, 0, 0) edge[line width=4, color=magenta] (3, 0, 0);

        \draw
            (0, 0, 1.8) edge[line width=4, color=red] (1, 0, 1)
            (1, 0, 1) edge[line width=4, color=blue] (2, 0, 1)
            (2, 0, 1) edge[line width=4, color=magenta] (3, 0, 1.8)
            (0, 0, 1.8) edge[line width=2] (0, 1, 1)
            (1, 0, 1) edge[line width=2] (0, 1, 1)
            (2, 0, 1) edge[line width=2] (0, 1, 1)
            (3, 0, 1.8) edge[line width=2] (0, 1, 1);

        \filldraw[black] (0, 1, 0) circle (1pt);
        \filldraw[black] (0, 0, 0) circle (1pt);
        \filldraw[black] (1, 0, 0) circle (1pt);
        \filldraw[black] (2, 0, 0) circle (1pt);
        \filldraw[black] (3, 0, 0) circle (1pt);

        \filldraw[black] (0, 1, 1) circle (1pt);
        \filldraw[black] (0, 0, 1.8) circle (1pt);
        \filldraw[black] (1, 0, 1) circle (1pt);
        \filldraw[black] (2, 0, 1) circle (1pt);
        \filldraw[black] (3, 0, 1.8) circle (1pt);

        \node[above left] at (0, 1, 0) {$z$};
        \node[below] at (1.5, 0, -0.1) {$P$};
        \node[above] at (0.1, 1, 1.05) {$(z, \omega)$};
    \end{tikzpicture}
\end{center}
            \caption*{(b) The graph of $\Psi$ lying over $\conv(z, P)$.}
        \end{minipage}
        \caption{}
        \label{triangles/figure/hyperplaneregulargraph}
    \end{figure}

    $\Psi$ is a continuous function $\conv(z, P) \longrightarrow \R$.
    As $\psi(\Z^{n} \cap P) \subseteq \Q$ and $\omega \in \Q$ it follows that $\Psi(\Z^n \cap \conv(z, P)) \subseteq \Q$.
    We are left to show the strict convexity of $\Psi$ with respect to $\cone(z, \mathcal S)$.

    We first show its convexity.
    For a cell $\sigma \in \mathcal S$, let $\psi_{\sigma}$ be the linear function $H \longrightarrow \R$ extending $\psi|_{\sigma}$.
    By convexity of $\psi$, we have
    \[
        \psi = \max_{\sigma} \psi_{\sigma}, \tag{$*$}
    \]
    where $\sigma$ ranges over the cells of $\mathcal S$.
    Define $\Psi_{\sigma}$ to be the linear extension of $\psi_{\sigma}$ with $\Psi_{\sigma}(z) = \omega$.
    Then we have
    \[
        \Psi|_{\conv(z, \sigma)} = \Psi_{\sigma}|_{\conv(z, \sigma)}.
    \]
    If we can show that
    \[
        \Psi = \max_{\sigma} \Psi_{\sigma}, \tag{$**$}
    \]
    where $\sigma$ ranges over the cells of $\mathcal S$, then we will have shown $\Psi$ is the max of a set of linear functions (which are convex), and hence that $\Psi$ is convex.

    It suffices to show that for $\sigma \in \mathcal S$ a cell and $y \in \conv(z, \sigma)$ we have $\Psi(y) = \Psi_{\sigma}(y) \geq \Psi_{\tau}(y)$ for any other cell $\tau$ of $\mathcal S$.
    Indeed, write $y = tx + (1 - t)z$ for $t \in [0, 1]$ and $x \in \sigma$.
    As $x \in \sigma$, $\psi(x) = \psi_{\sigma}(x)$.
    By $(*)$ we have $\psi_{\sigma}(x) \geq \psi_{\tau}(x)$.
    Then we compute
    \begin{align*}
        \Psi(y) &= \Psi_{\sigma}(y)\\
        &= t \psi_{\sigma}(x) + (1 - t) \omega\\
        &\geq t \psi_{\tau}(x) + (1 - t) \omega\\
        &= \Psi_{\tau}(y),
    \end{align*}
    and we have proven $(**)$ so convexity of $\Psi$ follows.

    We are left to show strict convexity with respect to $\cone(z, \mathcal S)$, i.e. that the domains of linearity of $\Psi$ are precisely the cells of $\cone(z, \mathcal S)$.
    These cells are of the form $\conv(z, \sigma)$ for $\sigma$ a cell of $\mathcal S$.
    $\Psi$ is by definition linear on these cells.
    Suppose then that there exist two cells $\conv(z, \sigma_1), \conv(z, \sigma_2)$ of $\cone(z, \mathcal S)$ so that $\Psi_{\sigma_1} = \Psi_{\sigma_2}$.
    As $\psi_{\sigma_i} = \Psi_{\sigma_i}|_{H}$, we have $\psi_{\sigma_1} = \psi_{\sigma_2}$.
    Therefore, because $\psi$ is strictly convex with respect to $\mathcal S$, we have that $\sigma_1 = \sigma_2$.
    Hence, $\conv(z, \sigma_1) = \conv(z, \sigma_2)$, so we have shown that $\Psi$ is strictly convex with respect to $\cone(z, \mathcal S)$, witnessing therefore the regularity of this subdivision.
\end{proof}

\begin{lemma}
\label{triangles/halfspaceregular}
    Let $H$ be a hyperplane in $\R^n$, which splits $\R^n$ into two closed half-spaces $H^{-}$ and $H^{+}$.
    Let $P$ be a polytope in $\R^n$ and define
    \begin{align*}
        P^{-} &= P \cap H^{-}\\
        P^{+} &= P \cap H^{+}\\
        P^{0} &= P \cap H.
    \end{align*}
    Suppose that $P^{0}$ is a facet of $P^{\pm}$ and that $P^{+} = \conv(z, P^{0})$ for some lattice point $z$ not in $H^{-}$.
    Suppose additionally that $P^{-}$ has a regular subdivision $\mathcal S^{-}$ which restricts to a regular subdivision $\mathcal S^{0}$ of $P^{0}$.
    Let $\mathcal S^{+} = \cone(z, \mathcal S^{0})$.
    Form the subdivision $\mathcal S$ of $P$ by gluing $\mathcal S^{+}$ and $\mathcal S^{-}$ along $\mathcal S^{0}$.
    We claim that $\mathcal S$ is regular.
    See \cref{triangles/figure/halfspaceregular}.
\end{lemma}
    \begin{figure}[thp]
        \begin{minipage}{0.5\textwidth}
            \centering
            \begin{center}
    \begin{tikzpicture}[scale=1.25]
        \fill[color=red, opacity=0.6] (0, 0) -- (1, 0) -- (1, -1) -- (0, -2) -- cycle {};
        \fill[color=blue, opacity=0.6] (1, 0) -- (1, -1) -- (2, 0) -- cycle {};

        \draw[line width=2] (0, 0) -- (1, 0);
        \draw[line width=2] (1, 0) -- (2, 0);

        \draw[line width=2] (0, 0) -- (0, -2);
        \draw[line width=2] (1, 0) -- (1, -1);
        \draw[line width=2] (2, 0) -- (0, -2);

        \draw (-1.5, 0) -- (0, 0);
        \draw (2, 0) -- (3.5, 0);

        \node[above] at (-1, 0) {$H$};
        \node at (-1, -1) {$H^{-}$};
    \end{tikzpicture}
\end{center}
            \caption*{(a) $P^{-} \subseteq H^{-}$ with regular subdivision $\mathcal S^{-}$.}
        \end{minipage}
        \begin{minipage}{0.5\textwidth}
            \centering
            \begin{center}
    \begin{tikzpicture}[scale=1.25]
        \fill[color=darkgreen, opacity=0.42] (0, 0) -- (1, 0) -- (0, 1) -- cycle {};
        \fill[color=pinkish, opacity=0.6] (2, 0) -- (1, 0) -- (0, 1) -- cycle {};

        \draw[line width=2] (0, 0) -- (1, 0);
        \draw[line width=2] (1, 0) -- (2, 0);

        \draw[line width=2] (0, 0) -- (0, 1);
        \draw[line width=2] (1, 0) -- (0, 1);
        \draw[line width=2] (2, 0) -- (0, 1);

        \draw (-1.5, 0) -- (0, 0);
        \draw (2, 0) -- (3.5, 0);

        \node[above] at (3, 0) {$H$};
        \node[above] at (0, 1.1) {$z$};
    \end{tikzpicture}
\end{center}
            \caption*{(b) $P^{+} = \conv(z, P^{0})$ with the subdivision $\mathcal S^{+} = \cone(z, \mathcal S^{0})$.}
        \end{minipage}

        \vspace{1.5cm}

        \begin{minipage}{\textwidth}
            \centering
            \begin{center}
    \begin{tikzpicture}[scale=1.25]
        \fill[color=red, opacity=0.6] (0, 0) -- (1, 0) -- (1, -1) -- (0, -2) -- cycle {};
        \fill[color=blue, opacity=0.6] (1, 0) -- (1, -1) -- (2, 0) -- cycle {};

        \fill[color=darkgreen, opacity=0.42] (0, 0) -- (1, 0) -- (0, 1) -- cycle {};
        \fill[color=pinkish, opacity=0.6] (2, 0) -- (1, 0) -- (0, 1) -- cycle {};

        \draw[line width=2] (0, 0) -- (1, 0);
        \draw[line width=2] (1, 0) -- (2, 0);

        \draw[line width=2] (0, 0) -- (0, -2);
        \draw[line width=2] (1, 0) -- (1, -1);
        \draw[line width=2] (2, 0) -- (0, -2);
        \draw[line width=2] (0, 0) -- (0, 1);
        \draw[line width=2] (1, 0) -- (0, 1);
        \draw[line width=2] (2, 0) -- (0, 1);

        \draw (-1.5, 0) -- (0, 0);
        \draw (2, 0) -- (3.5, 0);

        \node[above] at (1.5, 0.5) {$P^{+}$};
        \node at (1.5, -1) {$P^{-}$};
    \end{tikzpicture}
\end{center}
            \caption*{(c) The glued polytope $P = P^{-} \cup P^{+}$ with subdivision $\mathcal S = \mathcal S^{-} \cup \mathcal S^{+}$.}
        \end{minipage}
        \caption{}
        \label{triangles/figure/halfspaceregular}
    \end{figure}
\begin{proof}
    As in the proof of \cref{triangles/hyperplaneregular}, we let $\psi^{-}: P^{-} \longrightarrow \R$ witness the regularity of $\mathcal S$ and seek to extend $\psi^{-}$ linearly to some $\psi: P \longrightarrow \R$ which will witness the regularity of $\cone(z, \mathcal S)$.
    This time, we cannot choose $\psi(z) = \omega$ so freely, however.
    For each cell $\sigma$ of $\mathcal S^{-}$, let $\psi_{\sigma}$ be the linear function which agrees with $\psi^{-}$ on $\sigma$.
    Let $\omega$ be any rational number greater than $\psi_{\sigma}(z)$ for all cells $\sigma$ of $\mathcal S^{-}$.
    Define $\psi: P \longrightarrow \R$ to agree with $\psi^{-}$ on $P^{-}$, take the value $\omega$ on $z$, and linearly interpolate in between.
    See \cref{triangles/figure/halfspaceregularproof} for a visual of the way different choices of $\omega$ affect the convexity of $\psi$.

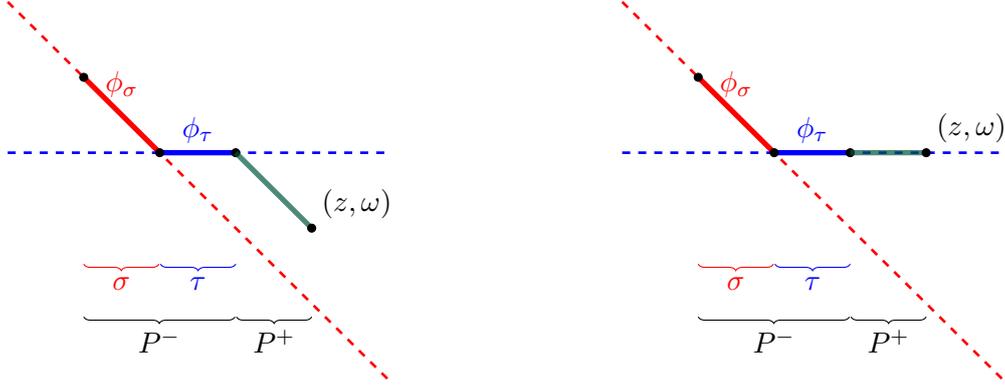
    \begin{figure}
        \begin{minipage}{0.5\textwidth}
            \centering
            \begin{center}
    \begin{tikzpicture}[scale=1]

        \draw[decoration={brace,mirror,raise=5pt},decorate] (0, -0.5) -- node[below=6pt] {$P^{-}$} (1.99, -0.5);
        \draw[decoration={brace,mirror,raise=5pt},decorate] (2.01, -0.5) -- node[below=6pt] {$P^{+}$} (3, -0.5);

        \draw[decoration={brace,mirror,raise=5pt},decorate, color=red] (0, 0.2) -- node[below=6pt] {$\sigma$} (0.99, 0.2);
        \draw[decoration={brace,mirror,raise=5pt},decorate, color=blue] (1.01, 0.2) -- node[below=6pt] {$\tau$} (2, 0.2);


        \draw[line width=2, color=red] (0, 2.5) -- (1, 1.5);
        \node[above, color=red] at (0.5, 2.1) {$\psi_{\sigma}$};
        \draw[line width=2, color=blue] (1, 1.5) -- (2, 1.5);
        \node[above, color=blue] at (1.5, 1.5) {$\psi_{\tau}$};

        \draw[line width=1, color=blue, dashed] (-1, 1.5) -- (4, 1.5);
        \draw[line width=1, color=red, dashed] (-1, 3.5) -- (4, -1.5);

        \filldraw[black] (0, 2.5) circle (1.5pt);
        \filldraw[black] (1, 1.5) circle (1.5pt);
        \filldraw[black] (2, 1.5) circle (1.5pt);

        \draw[line width=2, color=darkgreen, opacity=0.7] (2, 1.5) -- (3, 0.5);
        \filldraw[black] (3, 0.5) circle (1.5pt);
        \node[above right] at (3, 0.5) {$(z, \omega)$};
    \end{tikzpicture}
\end{center}
            \caption*{(a) The graph of $\psi$ if we choose $\omega$ strictly \\ below $\psi_\sigma(z)$ or $\psi_\tau(z)$.
            $\psi$ is not convex.}
        \end{minipage}
        \begin{minipage}{0.5\textwidth}
            \centering
            \begin{center}
    \begin{tikzpicture}[scale=1]

        \draw[decoration={brace,mirror,raise=5pt},decorate] (0, -0.5) -- node[below=6pt] {$P^{-}$} (1.99, -0.5);
        \draw[decoration={brace,mirror,raise=5pt},decorate] (2.01, -0.5) -- node[below=6pt] {$P^{+}$} (3, -0.5);

        \draw[decoration={brace,mirror,raise=5pt},decorate, color=red] (0, 0.2) -- node[below=6pt] {$\sigma$} (0.99, 0.2);
        \draw[decoration={brace,mirror,raise=5pt},decorate, color=blue] (1.01, 0.2) -- node[below=6pt] {$\tau$} (2, 0.2);


        \draw[line width=2, color=red] (0, 2.5) -- (1, 1.5);
        \node[above, color=red] at (0.5, 2.1) {$\psi_{\sigma}$};
        \draw[line width=2, color=blue] (1, 1.5) -- (2, 1.5);
        \node[above, color=blue] at (1.5, 1.5) {$\psi_{\tau}$};

        \draw[line width=1, color=blue, dashed] (-1, 1.5) -- (4, 1.5);
        \draw[line width=1, color=red, dashed] (-1, 3.5) -- (4, -1.5);

        \filldraw[black] (0, 2.5) circle (1.5pt);
        \filldraw[black] (1, 1.5) circle (1.5pt);
        \filldraw[black] (2, 1.5) circle (1.5pt);

        \draw[line width=2, color=darkgreen, opacity=0.7] (2, 1.5) -- (3, 1.5);
        \filldraw[black] (3, 1.5) circle (1.5pt);
        \node[above right] at (3, 1.5) {$(z, \omega)$}; 
    \end{tikzpicture}
\end{center}
            \caption*{(b) The graph of $\psi$ if we choose $\omega$ non-strictly above $\psi_\sigma(z)$ and $\psi_\tau(z)$. $\psi$ is non-strictly convex.}
        \end{minipage}
        
        \vspace{1cm}

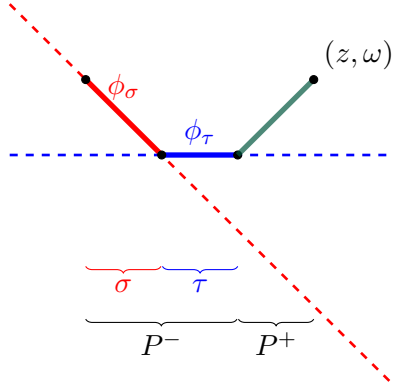
        \begin{minipage}{\textwidth}
            \centering
            \begin{center}
    \begin{tikzpicture}[scale=1]

        \draw[decoration={brace,mirror,raise=5pt},decorate] (0, -0.5) -- node[below=6pt] {$P^{-}$} (1.99, -0.5);
        \draw[decoration={brace,mirror,raise=5pt},decorate] (2.01, -0.5) -- node[below=6pt] {$P^{+}$} (3, -0.5);

        \draw[decoration={brace,mirror,raise=5pt},decorate, color=red] (0, 0.2) -- node[below=6pt] {$\sigma$} (0.99, 0.2);
        \draw[decoration={brace,mirror,raise=5pt},decorate, color=blue] (1.01, 0.2) -- node[below=6pt] {$\tau$} (2, 0.2);


        \draw[line width=2, color=red] (0, 2.5) -- (1, 1.5);
        \node[above, color=red] at (0.5, 2.1) {$\psi_{\sigma}$};
        \draw[line width=2, color=blue] (1, 1.5) -- (2, 1.5);
        \node[above, color=blue] at (1.5, 1.5) {$\psi_{\tau}$};

        \draw[line width=1, color=blue, dashed] (-1, 1.5) -- (4, 1.5);
        \draw[line width=1, color=red, dashed] (-1, 3.5) -- (4, -1.5);

        \filldraw[black] (0, 2.5) circle (1.5pt);
        \filldraw[black] (1, 1.5) circle (1.5pt);
        \filldraw[black] (2, 1.5) circle (1.5pt);

        \draw[line width=2, color=darkgreen, opacity=0.7] (2, 1.5) -- (3, 2.5);
        \filldraw[black] (3, 2.5) circle (1.5pt);
        \node[above right] at (3, 2.5) {$(z, \omega)$};
    \end{tikzpicture}
\end{center}
            \caption*{(c) The graph of $\psi$ if we choose $\omega$ strictly above $\psi_\sigma(z)$ and $\psi_\tau(z)$. $\psi$ is strictly convex.}
        \end{minipage}

        \caption{Various choices of $\omega$ in the proof of \cref{triangles/halfspaceregular} and the resulting convexity of $\psi$.}
        \label{triangles/figure/halfspaceregularproof}
    \end{figure}
                        
    $\psi$ is continuous, and as $\psi^{-}$ witnessed the regularity of $\mathcal S^{-}$, $\psi^{-}(\Z^{n} \cap P^{-}) \subseteq \Q$.
    So as $\omega \in \Q$, we have $\psi(\Z^{n} \cap P) \subseteq \Q$.
    We are left to show strict convexity with respect to $\mathcal S$.

    First, we show that the condition $\psi(z) > \psi_{\sigma}(z)$ for all cells $\sigma$ of $\mathcal S^{-}$ implies that $\psi(y) \geq \psi_{\sigma}(y)$ for any $y \in P^{+}$.
    Indeed, suppose $\psi(y) < \psi_{\sigma}(y)$ for some cell $\sigma \in \mathcal S^{-}$ and $y \in P^{+}$.
    We may assume $y \neq z$ and $y \notin P^{0}$.
    Then the line connecting $z$ and $y$ intersects the hyperplane $H$ at a unique point $x$.
    We have the following inequalities
    \begin{align*}
        (\psi - \psi_{\sigma})(x) &\geq 0\\
        (\psi - \psi_{\sigma})(y) &< 0\\
        (\psi - \psi_{\sigma})(z) &> 0
    \end{align*}
    where $x, y, z$ are distinct and collinear.
    Both $\psi_{\sigma}$ and $\psi$ are linear on the line segment $\ell$ from $x$ to $y$ to $z$.
    The linear function $(\psi - \psi_{\sigma})|_{\ell}$ cannot simultaneously satisfy all three of these inequalities, so $\psi(y) \geq \psi_{\sigma}(y)$ for all points $y \in P^{+}$ and cells $\sigma \in \mathcal S^{-}$.

    We now show convexity itself.
    Let $\psi^+ = \psi|_{P^+}$.
    $\psi^{-}$ is convex by hypothesis, and $\psi^{+}$ is convex by \cref{triangles/hyperplaneregular}.
    Hence, we need only check convexity for line segments connecting points $x^{-} \in P^{-}$ and $x^{+} \in P^{+}$ which are not in $P^{0}$.
    Let $\ell: [-1, 1] \longrightarrow P$ be the line segment so that $\ell(-1) = x^{-}$ and $\ell(1) = x^{+}$.
    We seek to show $\psi \circ \ell$ has increasing slopes.
    Let $t_0 \in (-1, 1)$ be the unique value where $\ell(t_0) \in P^0$.
    By convexity of $\psi^{\pm}$, $\psi \circ \ell$ has increasing slopes on $[-1, t_0]$ and on $[t_0, 1]$.
    We claim that the slope of $\psi \circ \ell$ at $t_0 - \varepsilon$ is less than the slope of $\psi \circ \ell$ at $t_0 + \varepsilon$ for small $\varepsilon$.
    For all sufficiently small $\varepsilon$, $\ell(t_0 - \varepsilon)$ is in some fixed cell $\sigma \in \mathcal S^{-}$, so $\psi(\ell(t)) = \psi_{\sigma}(\ell(t))$ for $t \in [t_0 - \varepsilon, t_0]$.
    As discussed above, we know that $\psi(y) \geq \psi_\sigma(y)$ for any point $y \in P^{+}$ and cell $\sigma \in \mathcal S^{-}$.
    Hence, $\psi(\ell(t)) \geq \psi_{\sigma}(\ell(t))$ for $t \in [t_0, t_0 + \varepsilon]$.
    As such, $\psi \circ \ell$ has larger slope at $t_0 + \varepsilon$ than $t_0 - \varepsilon$, proving convexity.

    Finally, we show the convexity is strict with respect to $\mathcal S$.
    Certainly, $\psi$ is linear on the cells of $\mathcal S$.
    For cells $\sigma_1$, $\sigma_2$ in $\mathcal S$, let $\psi_{\sigma_i}$ denote the linear function which agrees with $\psi$ on $\sigma_i$.
    Suppose $\psi_{\sigma_1} = \psi_{\sigma_2}$.
    If both $\sigma_i$ are in $P^{-}$ then by regularity of $\mathcal S^{-}$, $\sigma_1 = \sigma_2$.
    If $\sigma_1$ is in $P^{-}$ and $\sigma_2$ is in $P^{+}$, then we have by definition that $\psi_{\sigma_2}(z) = \psi(z) > \psi_{\sigma_1}(z)$, contradicting the hypothesis that $\psi_{\sigma_1} = \psi_{\sigma_2}$
    Finally, suppose both $\sigma_i$ are in $P^{+}$.
    Then $\sigma_i = \conv(z, \tau_i)$ for unique cells $\tau_i$ of $\mathcal S^{0}$.
    As $\mathcal S^0$ is regular and $\psi_{\sigma_1}|_{H} = \psi_{\sigma_2}|_{H}$, it follows that $\tau_1 = \tau_2$.
    Hence, $\sigma_1 = \sigma_2$ and we conclude strict convexity.
\end{proof}

Now, consider a regular, unimodular, star-shaped triangulation $\mathcal T$ of $P_2^{(n)}$.
By \cref{triangles/hyperplaneregular}, $\cone\parens{e_n^{(n + 1)}, \mathcal T}$ is a regular triangulation of $\conv\parens{e_n^{(n + 1)}, P_2^{(n)}}$.
By \cref{triangles/halfspaceregular}, gluing this to $\conv\parens{w_1^{(n + 1)}, P_2^{(n)}}$ yields a regular triangulation $\mathcal T'$ of $\conv\parens{w_1^{(n + 1)}, \conv\parens{e_n^{(n + 1)}, P_2^{(n)}}}$, which is exactly $P_1^{(n + 1)}$.
By \cref{triangles/coneunimodular}, $\mathcal T'$ is unimodular.
Furthermore, star-shapedness is preserved by taking cones.
We have shown therefore that to find a regular, unimodular, star-shaped triangulation of $P_1^{(n + 1)}$ it suffices to find a regular, unimodular, star-shaped triangulation of $P_2^{(n)}$.
That is, \cref{triangles/triangulating1} is reduced to the following.
\begin{proposition}
\label{triangles/triangulating2}
    For all $n$, there is a regular, unimodular, star-shaped triangulation of $P_2^{(n)}$.
\end{proposition}

In the following section, we will construct these triangulations.

\subsection{Triangulating $P_2$}
\label{section/triangles/triangulatingP2}

We seek a recursive argument to find regular, unimodular, star-shaped triangulations of the family $P_2$.
Observe that the family of polytopes $P_2$ all sit inside each other as hyperplane sections setting the last coordinate equal to 0.
That is,
\[
    \{x_n = 0\} \cap P_2^{{(n + 1)}} = P_2^{{(n)}} \times \{0\},
\]
which can be argued similarly to \cref{triangles/2inside1}.
Unlike the similar relationship between $P_1$ and $P_2$ in \cref{triangles/2inside1}, where precisely two new lattice points are added, $P_2^{(n + 1)}$ has many more lattice points than $P_2^{(n)}$.
As such, a recursive argument using this containment is inconvenient.
Instead, we note that $P_2$ is (remarkably) isomorphic to its \hyperlink{background/polardual}{polar dual} as a lattice polytope, and that the dual polytopes satisfy a much cleaner recurrence.

\begin{lemma}
\label{triangles/selfduality}
    For all $n$, the linear map $T: P_2^{(n)} \longrightarrow \widecheck{P}_{2}^{(n)}$ given by
    \[
        e_i \mapsto (-1, \dots, s_i - 1, \dots, -1),
    \]
    where $s_i - 1$ is in the $i^{th}$ position (indexed from $0$), is an isomorphism of lattice polytopes.
\end{lemma}
\begin{proof}
    We first ensure that $T$ actually sends $P_2^{(n)}$ to $\widecheck{P}_2^{(n)}$.
    Recall that $P_2^{(n)}$ has vertices $e_0, \dots, e_{n - 1}, w_2$.
    We claim then that $T(e_0), \dots, T(e_{n - 1}), T(w_2)$ are the vertices of $\widecheck{P}_2^{(n)}$.
    By polar duality (\cref{background/polarduality}), this means we must show that these vertices define the half-spaces cutting out $P_2^{(n)}$.

    We first compute $T(w_2) = (-1, \dots, -1)$.
    Indeed, by linearity of $T$, this means we must show
    \[
        -1 = -\frac{d_2}{s_{i}}(s_{i} - 1) + \sum_{j \neq i} \frac{d_2}{s_i}
    \]
    for all $0 \leq i \leq n - 1$, where $d_2 = d_2^{(n)} = s_{n} - 1$.
    Sylvester's sequence satisfies the following recurrence:
    \[
	\sum_{i=0}^{n-1} \frac{d_2}{s_i} = d_2 - 1.
    \]
    Hence, for all $i$,
    \begin{align*}
        -\frac{d_2}{s_i} (s_i - 1) + \sum_{j \neq i} \frac{d_2}{s_i} &= (d - 1) - \parens{\frac{d}{s_{i}} + \frac{d}{s_{i}}(s_{i} - 1)}\\
        &= -1,
    \end{align*}
    as desired.

    We now seek to show that $T(e_0), \dots, T(e_{n - 1}), T(w_2)$ form the half-spaces cutting out $P_2^{(n)}$.
    Indeed, to show a simplex $\conv\parens{v_0, \dots, v_n}$ is cut out by the $n + 1$ half-spaces $\{\Angle{x, u_i} + 1 \geq 0\}$, we must show that each $v_i$ satisfies equality in $n$ of these inequalities and satisfies the strict inequality in the last.
    In our context, $v_0 = e_0, \dots, v_{n - 1} = e_{n - 1}, v_n = w_2$, and $u_i = T(v_i)$.
    
    First, consider $v_i$ for $i < n$.
    Then $\Angle{v_i, u_j}$ is the $i^{th}$ component of $u_j$, where we index components starting from $0$.
    For $i \neq j$, the $i^{th}$ component of $u_j$ is $-1$, hence we have equality $\Angle{v_i, u_j} + 1 = 0$.
    For $i = j$, the $i^{th}$ component of $u_i$ is $s_i - 1$, so $\Angle{v_i, u_i} + 1 > 0$.

    Now, consider $v_n$.
    For $u_n$, we have
    \begin{align*}
        \Angle{v_n, u_n} &= \sum_{i = 0}^{n - 1} \frac{d_2}{s_i}\\
        &>0.
    \end{align*}
    For $u_j$ with $j < n$, the above computation that $T(v_n) = (-1, \dots, -1)$ showed exactly that $\Angle{v_n, u_j} + 1 = 0$.

    We have computed the vertices of $\widecheck{P}_2^{(n)}$ to be precisely the image of the vertices of $P_2^{(n)}$ under $T$.
    To conclude that these simplices are isomorphic as lattice polytopes, we must show that $T$ has determinant $\pm 1$.
    The matrix representation of $T$ is
    \[
        \begin{pmatrix}
            1 & -1 & \dots & -1\\
            -1 & 2 & & -1\\
            \vdots & & \ddots & \vdots\\
            -1 & -1 & \dots & s_{n - 1} - 1
        \end{pmatrix}.
    \]
    We compute this determinant by simultaneous row and column operations to clear the first row and column.
    For instance, doing this once yields:

    \begin{align*}
            \det(T)
        &=
            \det
            \begin{pmatrix}
                1 & 0 & 0 & \dots & 0\\
                0 & 1 & -2 & \dots & -2\\
                0 & -2 & 3 & & -1\\
                \vdots & \vdots & & \ddots & \vdots\\
                0 & -2 & -1 & \dots & s_{n - 1} - 2
            \end{pmatrix}
        \\
        &=
            \det
            \begin{pmatrix}
                1 & -2 & \dots & -2\\
                -2 & 3 & & -1\\
                \vdots & & \ddots & \vdots\\
                -2 & -1 & \dots & s_{n - 1} - 2
            \end{pmatrix}
            ,
    \end{align*}
    where we are positioned again to clear the first row and column.
    Performing this process $k$ times from $T$ yields
    \[
        \det(T) =
        \det
        \begin{pmatrix}
            1 & 1 - s_k & \dots & 1 - s_k\\
            1 - s_k & s_{k + 1} + (1 - s_k) & & 1 - s_k\\
            \vdots & & \ddots & \vdots\\
            1 - s_k & 1 - s_k & \dots & s_{n - 1} + (1 - s_k)
        \end{pmatrix}
        ,
    \]
    which can be verified inductively using the following relation of Sylvester's sequence:
    \[
        (1 - s_k) - (1 - s_k)^2 = 1 - s_{k + 1}.
    \]
    Taking $k = n - 1$ we therefore determine $\det(T) = 1$, proving that $T$ yields an isomorphism of lattice polytopes $P_2^{(n)} \xlongrightarrow{\sim} \widecheck{P}_2^{(n)}$.
\end{proof}

We have therefore reduced the proof of \cref{triangles/triangulating2} (and hence of \cref{triangles/triangulating1}) to the following.
\begin{proposition}
\label{triangles/triangulating2dual}
    For all $n$, there is a regular, unimodular, star-shaped triangulation $\mathcal T^{(n)}$ of $\widecheck{P}_2^{(n)}$.
\end{proposition}
We will proceed to construct these triangulations through the remainder of this section.
This is notably simpler than working directly with $P_2$ because the family $\widecheck{P}_2$ has a much cleaner recursive structure.
First, its faces are simple to describe.
As $e_0, \dots, e_n$ are vertices of $P_2^{(n + 1)}$, it follows that the coordinate hyperplanes $\{x_i = -1\}$ define all but one of the faces of $\widecheck{P}_2^{(n + 1)}$.
The final ``slanted" face is defined by $\curly{\sum_{i = 0}^{n} -\frac{d_2^{(n + 1)}}{s_i} x_i = -1}$.
Furthermore, the map $\R^n \longrightarrow \R^{n + 1}$ taking $x \mapsto (x, -1)$ embeds $\widecheck{P}_2^{(n)}$ isomorphically as the face $\{x_n = -1\} \cap \widecheck{P}_2^{(n + 1)}$.
We refer to this as the ``bottom" face of $\widecheck{P}_2^{(n + 1)}$.
See \cref{triangles/figure/polarrecursive} for a visual of this relationship.
\begin{figure}
    \centering
    \usetikzlibrary{quotes}

\begin{tikzpicture}[
        scale=2,
        every edge quotes/.style = {auto, sloped}
    ]
    \draw
        (-1, -1) edge[line width=2, auto=right, "$\widecheck{P}_2^{(n)} \times \{-1\}$"] (1, -1)
        (-1, -1) edge[line width=2] (-1, 2)
        (-1, 2) edge[line width=2] (1, -1)
        (-2, -1) edge (2, -1);

    \node[above] at (2, -1) {$\{x_n = -1\}$};
    \node[above right] at (0, 0.5) {$\sum -\frac{d_2^{(n + 1)}}{s_i} x_i = -1$};
    \node[above] at (-1, 2) {$(-1, \dots, -1, s_n - 1)$};
\end{tikzpicture}
    \caption{A schematic representation of $\widecheck{P}_2^{(n + 1)}$ and its recursive face structure.
	Literally speaking, the triangle depicted is $\widecheck{P}_2^{(2)}$.}
    \label{triangles/figure/polarrecursive}
\end{figure}
This follows directly from the computation in the proof of \cref{triangles/selfduality} of the vertices of $\widecheck{P}_2^{(n + 1)}$.
Namely, we showed that
\[
    \widecheck{P}_2^{(n + 1)} = \conv\parens{
        \begin{pmatrix}
            -1\\
            -1\\
            \vdots\\
            -1
        \end{pmatrix}
        ,
        \begin{pmatrix}
            1\\
            -1\\
            \vdots\\
            -1
        \end{pmatrix}
        ,
        \begin{pmatrix}
            -1\\
            2\\
            \vdots\\
            -1
        \end{pmatrix}
        ,\dots, 
        \begin{pmatrix}
            -1\\
            -1\\
            \vdots\\
            s_{n} - 1
        \end{pmatrix}
    }
    ,
\]
and setting $x_n = -1$ yields $\widecheck{P}_2^{(n)}$.

An immediate benefit of the recursive structure of $\widecheck{P}_2$ is that we can understand the structure of the lattice points within these simplices systematically.
\begin{lemma}
\label{triangles/columns}
    Let $\pi: \R^{n + 1} \longrightarrow \R^n$ be projection onto the first $n$ coordinates.
    For every $x \in \widecheck{P}_2^{(n + 1)}$, $\pi(x) \in \widecheck{P}_2^{(n)}$.
\end{lemma}
\begin{proof}
    Let $x \in \widecheck{P}_2^{(n + 1)}$.
    Then we may write $x$ as a convex combination $x = \sum \lambda_i u_i$, where the $u_i$ are the vertices of $\widecheck{P}_2^{(n + 1)}$ computed above.
    Take $u_n = (-1, \dots, -1, s_{n} - 1)$ and $u_0 = (-1, \dots, -1)$.
    Then $\pi(u_n) = \pi(u_0), \dots, \pi(u_{n - 1})$ are the vertices of $\widecheck{P}_2^{(n)}$.
    Hence, $\pi(x) = \sum \lambda_i \pi(u_i)$ is a convex combination of vertices of $\widecheck{P}_2^{(n)}$.
    So $\pi(x)$ is in $\widecheck{P}_2^{(n)}$, as desired.
\end{proof}

For $x \in \widecheck{P}_2^{(n + 1)}$, we say $x$ lies ``above" $(\pi(x), -1)$.
From \cref{triangles/columns}, lattice points of $\widecheck{P}_2^{(n + 1)}$ come in columns above the lattice points of $\widecheck{P}_2^{(n)} \times \{-1\}$.
We determine now the highest point on each such column, ordered by $x_n$-coordinate (recall $\widecheck{P}_2^{(n + 1)}$ is a polytope in $\R^{n + 1}$ with coordinates $(x_0, \dots, x_n)$).
\begin{lemma}
\label{triangles/highestpoints}
    Let $y \in \widecheck{P}_2^{(n)}$.
    Let $x$ be the lattice point of $\widecheck{P}_2^{(n + 1)}$ satifying $\pi(x) = y$ with $x_n$ maximal.
    \begin{enumerate}[label = (\alph*)]
        \item If $y = (-1, \dots, -1)$ then $x = (-1, \dots, -1, s_n - 1)$.
        \item Otherwise, $x$ satisfies
        \[
            x_n = -\sum_{i = 0}^{n - 1} \frac{s_n - 1}{s_i} x_i.
        \]
    \end{enumerate}
\end{lemma}
See \cref{triangles/figure/highestpoints}.

\begin{figure}
    \centering
    \usetikzlibrary{quotes}

\begin{tikzpicture}[
        scale=2,
        every edge quotes/.style = {auto, sloped}
    ]
    \draw
        (-1, -1) edge[line width=2] (1, -1)
        (-1, -1) edge[line width=2] (-1, 2)
        (-1, 2) edge[line width=2] (1, -1)
        (-1, 1) edge[line width=2, auto=right, "$H = \curly{x_n = -\sum \frac{s_n - 1}{s_i} x_i}$"] (1, -1);

    \filldraw[black] (-1, -1) circle (1pt);
    \filldraw[black] (-1, 0) circle (1pt);
    \filldraw[black] (-1, 1) circle (1pt);
    \filldraw[black] (-1, 2) circle (1pt);
    \filldraw[black] (0, -1) circle (1pt);
    \filldraw[black] (0, 0) circle (1pt);
    \filldraw[black] (1, -1) circle (1pt);

    \draw[decoration={brace,raise=5pt},decorate] (-1, -1) -- node[left=6pt] {$P^{\leq H}$} (-1, 0.99);
    \draw[decoration={brace,raise=5pt},decorate] (-1, 1.01) -- node[left=6pt] {$P^{\geq H}$} (-1, 2);
    \draw[decoration={brace,raise=15pt, mirror},decorate] (-1, -1) -- node[below=21pt] {$\widecheck{P}_2^{(n)}$} (1, -1);

    \node[above] at (-1, 2) {$(-1, \dots, -1, s_n - 1)$};
    \node[below] at (-1, -1) {$(-1, \dots, -1)$};
\end{tikzpicture}
    \caption{A schematic representation of $\widecheck{P}_2^{(n + 1)}$ showing its lattice points and how it is cut by the hyperplane $H$.
	Literally speaking, the triangle depicted is $\widecheck{P}_2^{(2)}$.}
    \label{triangles/figure/highestpoints}
\end{figure}
\begin{proof}
    \begin{enumerate}[label = (\alph*)]
        \item Any $z \in \widecheck{P}_2^{(n + 1)}$ is of the form $\sum \lambda_i u_i$ where $\lambda_i \geq 0$, $\sum \lambda_i = 1$, and $u_0, \dots, u_{n + 1}$ are the vertices of $\widecheck{P}_2^{(n + 1)}$ computed in the course of proving \cref{triangles/selfduality}.
        Take $u_{n + 1} = (-1, \dots, -1, s_n - 1)$ so that the other $u_i$ have last component $-1$.
        Hence, the last component of $z$ is $\lambda_n (s_n - 1) - \sum_{i = 0}^{n - 1} \lambda_i$, which is at most $s_n - 1$.
        So in fact $x = u_{n + 1}$ is the highest point on the entirety of $\widecheck{P}_2^{(n + 1)}$, and is thus the highest point in its column.

        \item The equation given specifies the last coordinate in terms of the previous ones.
        Indeed, write
        \[
            x = \parens{y_0, \dots, y_{n - 1}, -\sum_{i = 0}^{n - 1} \frac{s_n - 1}{s_i} y_i}.
        \]
        We seek to show that $x$ is indeed the highest point above $y$.
        
        First, we show that $x \in \widecheck{P}_2^{(n + 1)}$ by checking that all $x_i \geq -1$ and $ - \sum \frac{s_{n + 1} - 1}{s_i} x_i \geq -1$.
        
        Let's begin with the coordinate inequalities.
        For $i < n$, we defined $x_i = y_i$, and we know $y_i \geq -1$ as $y \in \widecheck{P}_2^{(n)}$.
        Now, as $y \in \widecheck{P}_2^{(n)}$ we have
        \[
            -\sum_{i = 0}^{n - 1} \frac{s_n - 1}{s_i} y_i \geq -1.
        \]
        The left hand side is precisely $x_n$, so $x_n \geq -1$.

        Now, we show the last inequality
        \[
            -\sum_{i = 0}^{n} \frac{s_{n + 1} - 1}{s_i} x_i \geq -1.
        \]
        Indeed,
        \begin{align*}
            -\sum_{i = 0}^{n} \frac{s_{n + 1} - 1}{s_i} x_i &= s_n \sum_{i = 0}^{n - 1} -\frac{s_n - 1}{s_i} x_i - (s_n - 1) x_n\\
            &= s_n x_n - (s_n - 1) x_n\\
            &= x_n \tag{$*$}\\
            &\geq -1,
        \end{align*}
        so we have shown $x \in \widecheck{P}_2^{(n + 1)}$.

        Finally, we show that nothing above $x$ is in $\widecheck{P}_2^{(n + 1)}$.
        Indeed, consider $x' = x + e_n$.
        By assumption, one of the $x_i$ for $i < n$ is strictly greater than $-1$.
        Hence, 
        \begin{align*}
            -\sum_{i = 0}^{n - 1} \frac{s_n - 1}{s_i} x_i &< \sum_{i = 0}^{n - 1} \frac{s_n - 1}{s_i}.
        \end{align*}
        The left side of this inequality is $x_n$.
        The right side can be computed as
        \begin{align*}
            \sum_{i = 0}^{n - 1} \frac{s_n - 1}{s_i} &= (s_n - 1)\parens{1 - \frac{1}{s_n - 1}}\\
            &= s_n - 2.
        \end{align*}
        So we conclude that $x_n < s_n - 2$.

        Now to show $x' \notin \widecheck{P}_2^{(n + 1)}$ we will show
        \[
            -\sum_{i = 0}^{n} \frac{s_{n + 1} - 1}{s_i} x_i' < -1.
        \]
        Using $(*)$, the left hand side may be computed as
        \begin{align*}
            -\sum_{i = 0}^{n} \frac{s_{n + 1} - 1}{s_i} x_i' &= x_n - \frac{s_{n + 1} - 1}{s_n}\\
            &= x_n - (s_n - 1)\\
            &< (s_n - 2) - (s_n - 1)\\
            &= -1.
        \end{align*}
        We have shown $x' \notin \widecheck{P}_2^{(n + 1)}$, so $x$ is indeed the highest element in its column.
    \end{enumerate}
\end{proof}

With this precise recursive relationship between the lattice points of $\widecheck{P}_2^{(n)}$ and $\widecheck{P}_2^{(n + 1)}$, we seek a recursive algorithm to find a regular, unimodular, star-shaped triangulation of $\widecheck{P}_2^{(n + 1)}$ given such a triangulation of $\widecheck{P}_2^{(n)}$.
To do so, we use the notion of the \emph{pulling refinement} due to Haase, Paffenholz, Piechnik, and Santos in \cite{unimodular}.

\begin{definition}
    \hypertarget{triangles/pulling}
    Let $\mathcal S$ be a subdivision of a polytope $P \subseteq \R^d$ and $m \in P \cap \Z^d$.
    The \emph{pulling refinement} $\mathrm{pull}_m(\mathcal S)$ is defined by replacing every $F \in \mathcal S$ containing $m$ by $\conv(m, F')$ for every face $F' \leq F$ which does \emph{not} contain $m$.
\end{definition}

In \cite{unimodular}, the authors prove the following two key properties of the pulling refinement, which we restate here for convenience.

\begin{lemma}[Lemma 2.1 in \cite{unimodular}]
\label{triangles/pullingregular}
    If $\mathcal S$ is a regular subdivision of a polytope $P$ and $m \in P$ a lattice point, then $\mathrm{pull}_m(\mathcal S)$ is also a regular subdivision.
\end{lemma}

The proof of this proceeds by taking a witness $\psi$ to the regularity of $\mathcal S$ and producing a witness $\psi'$ to the regularity of $\mathrm{pull}_{m}(\mathcal S)$.
$\psi'$ is the piecewise linear function which agrees with $\psi$ on every lattice point of $P$ except $m$, and has $\psi'(m) = \psi(m) - \varepsilon$ for some sufficiently small $\varepsilon > 0$ (this, incidentally, explains why this is referred to as ``pulling" $m$ -- we take the graph of $\psi$ and pull $(m, \psi(m))$ down by $\varepsilon$).
To ensure $\Q$-linearity, we can simply take $\varepsilon$ to be of the form $1/N$ for some sufficiently large $N$.

\begin{definition}
    Let $T: \R^n \longrightarrow \R^m$ be linear and let $P$ be a polytope in $\R^n$.
    Set $Q = T(P)$.
    For a polytope $\sigma \subseteq Q$, we define its \emph{pullback} to be
    \[
        T^*(\sigma) = T^{-1}(\sigma) \cap P,
    \]
    and for a subdivision $\mathcal S$ of $Q$, we define the \emph{pullback subdivision} to be
    \[
        T^* \mathcal S = \{T^{*}(\sigma) \mid \sigma \in \mathcal S\}.
    \]
\end{definition}

\begin{proposition}[Theorem 2.8 in \cite{unimodular}]
\label{triangles/pullingunimodular}
    Let $\pi: \R^{n + 1} \longrightarrow \R^{n}$ be a projection so that we have $\pi(\Z^{n + 1}) = \Z^{n}$.
    Let $P$ be a lattice polytope in $\R^{n + 1}$ and set $Q = \pi(P)$.
    Suppose $\mathcal T$ is a unimodular triangulation of $Q$ such that $\pi^* \mathcal T$ is a \emph{lattice} subdivision.
    Let $\mathcal T'$ be a refinement of $\pi^{*} \mathcal T$ arising by pulling $\mathcal S$ at all lattice points in $P$ in any order.
    Then $\mathcal T'$ is a unimodular triangulation of $P$.
\end{proposition}

Combining these, we may use pulling refinements to lift regular, unimodular triangulations along projections.
(We will see in \cref{triangles/starshapedbelow} and the discussion afterwards how to preserve the \hyperlink{background/starshaped}{star-shaped} condition.)
We are tempted then to simply apply this to $P = \widecheck{P}_2^{(n + 1)}$ and $Q = \widecheck{P}_2^{(n)}$.
However, given a unimodular triangulation $\mathcal T$ of $\widecheck{P}_{2}^{(n)}$ and the map $\pi: \R^{n + 1} \longrightarrow \R^n$ omitting the last coordinate, $\pi^* \mathcal T$ is not a \emph{lattice} subdivision.
Indeed, the topmost face of $\widecheck{P}_2^{(n + 1)}$, defined by $\curly{\sum -\frac{d_2}{s_i} x_i = -1}$, has no lattice points other than its vertices.
This can be verified by the formula for the highest lattice points per column given in \cref{triangles/highestpoints}.

Nevertheless, this method does get us quite close to the top.
Consider the hyperplane $H = \curly{x_n = -\sum \frac{s_n - 1}{s_i}x_i} \subseteq \R^{n + 1}$ introduced in \cref{triangles/highestpoints} and let
\begin{align*}
    P^{= H} &= \curly{x \in \widecheck{P}_2^{(n + 1)} \bigg| \; x_n = -\sum_{i = 0}^{n} \frac{s_n - 1}{s_i} x_i}\\
    P^{\leq H} &= \curly{x \in \widecheck{P}_2^{(n + 1)} \bigg| \; x_n \leq -\sum_{i = 0}^{n} \frac{s_n - 1}{s_i} x_i}\\
    P^{\geq H} &= \curly{x \in \widecheck{P}_2^{(n + 1)} \bigg| \; x_n \geq -\sum_{i = 0}^{n} \frac{s_n - 1}{s_i} x_i}.
\end{align*}
See \cref{triangles/figure/highestpoints}.

$\pi$ sends $P^{=H}$ isomorphically as a lattice polytope to $\widecheck{P}_2^{(n)}$ and $\pi(P^{\leq H}) = \widecheck{P}_2^{(n)}$.
Now, we see that a lattice triangulation $\mathcal T$ of $\widecheck{P}_2^{(n)}$ induces a \emph{lattice} subdivision $\pi^* \mathcal T$ of $P^{\leq H}$, so \cref{triangles/pullingunimodular} applies.
                                                                                            
The only lattice point in $P^{\geq H}$ which is not on $H$ is $z = (-1, \dots, -1, s_n - 1)$.
And furthermore, $P^{\geq H} = \conv(z, P^{= H})$.
Indeed, $P$ is convex and the condition $x_n \geq -\sum_{i = 0}^{n} \frac{s_n - 1}{s_i} x_i$ is preserved by convex combinations, so $P^{\geq H}$ is convex and hence contains $\conv(z, P^{= H})$.
For any $y \in P^{\geq H}$ other than $z$, consider the line connecting $z$ to $y$.
This line will intersect $P^{= H}$ at a unique point $x$.
$y$ is therefore a convex combination of $x \in P^{= H}$ and $z$.
We have thus proven the equality $P^{\geq H} = \conv(z, P^{= H})$, so we may triangulate $P^{\geq H}$ via $\cone\parens{z, \parens{\pi^*\mathcal T}|_H}$.
\begin{lemma}
\label{triangles/aboveunimodular}
    If $\mathcal T$ is unimodular, the subdivision $\cone(z, (\pi^* \mathcal T)|_H)$ of $P^{\geq H}$ is unimodular.
\end{lemma}
\begin{proof}
    The linear transformation
    \[
        T(x_0, \dots, x_n) = \parens{x_0, \dots, x_{n - 1}, x_n + \sum_{i = 0}^{n} \frac{s_n - 1}{s_i} x_i}
    \]
    has determinant $1$ and hence preserves unimodularity.
    For $x \in H$, $T(x) = (\pi(x), 0)$, so $T$ sends $P^{= H}$ isomorphically to $\widecheck{P}_2^{(n)} \times \{0\}$ and $(\pi^* \mathcal T)|_H$ isomorphically to $\mathcal T \times \{0\}$.
    Furthermore, $T(-1, \dots, -1, s_n - 1) = (-1, \dots, -1, 1)$.
    We deduce therefore that
    \[
	T\parens{P^{\geq H}} = \conv\parens{(-1, \dots, -1, 1), \widecheck{P}_2^{(n)} \times \{0\}}.
    \]
    Hence, $T$ sends $\cone(z, (\pi^* \mathcal T)|_H)$ isomorphically to $\cone((-1, \dots, -1, 1), \mathcal T \times \{0\})$.
    This is unimodular by \cref{triangles/conevolume}, as $\mathcal T$ is unimodular.
\end{proof}

We have thus far subdivided $P^{\leq H}$ via $\pi^* \mathcal T$ and $P^{\geq H}$ via $\conv(z, (\pi^* \mathcal T)|_H)$.
Gluing these together yields a subdivision $\mathcal S$ of $\widecheck{P}_2^{(n + 1)} = P^{\leq H} \cup P^{\geq H}$.
See \cref{triangles/figure/pullbackcone}.

\begin{lemma}
\label{triangles/gluedregular}
    If $\mathcal T$ is regular then the subdivision $\mathcal S$ obtained by gluing $\pi^* \mathcal T$ and $\conv(z, \pi^* \mathcal T|_H)$ along $H$ is regular.
\end{lemma}
\begin{proof}
    First, we show that $\pi^* \mathcal T$ is regular.
    Indeed, let $\psi: \widecheck{P}_2^{(n)} \longrightarrow \R$ witness regularity of $\mathcal T$.
    That is, $\psi$ is a continuous, piecewise $\Q$-linear function which is strictly convex with respect to $\mathcal T$.
    Then $\psi \circ \pi: P^{\leq H} \longrightarrow \R$ is continuous, piecewise $\Q$-linear, and strictly convex with respect to $\pi^* \mathcal T$, and thus witnesses its regularity.
    
    Now we are in the situation of \cref{triangles/halfspaceregular}, letting $P^{-} = P^{\leq H}$ subdivided via $\pi^{*} \mathcal T$ and $z = (-1, \dots, -1, s_{n} - 1)$.
    We therefore conclude that $\mathcal S$ is regular.
\end{proof}

So far, we have shown that given a regular, unimodular, star-shaped triangulation $\mathcal T^{(n)}$ of $\widecheck{P}_2^{(n)}$ we can form a regular subdivision $\mathcal S$ of $\widecheck{P}_2^{(n + 1)}$ whose restriction $\mathcal S|_{P^{\geq H}}$ is a unimodular triangulation.
See \cref{triangles/figure/pullbackcone}.

\begin{figure}
    \centering
    \usetikzlibrary{quotes}

\begin{tikzpicture}[
        scale=2,
    ]
    \fill[color=red, opacity=0.6] (-1, -1) -- (0, -1) -- (0, 0) -- (-1, 1) -- cycle {};
    \fill[color=blue, opacity=0.6] (0, -1) -- (0, 0) -- (1, -1) -- cycle {};
    \draw
        (-1, -1) edge[line width=2] (1, -1)
        (-1, -1) edge[line width=2] (-1, 1)
        (0, -1) edge[line width=2] (0, 0)
        (-1, 1) edge[line width=2] (1, -1);

    \filldraw[black] (-1, -1) circle (1pt);
    \filldraw[black] (-1, 0) circle (1pt);
    \filldraw[black] (-1, 1) circle (1pt);
    \filldraw[black] (0, -1) circle (1pt);
    \filldraw[black] (0, 0) circle (1pt);
    \filldraw[black] (1, -1) circle (1pt);

    \fill[color=darkgreen, opacity=0.42] (-1, 1.2) -- (0, 0.2) -- (-1, 2.2) -- cycle {};
    \fill[color=pinkish, opacity=0.6] (1, -0.8) -- (0, 0.2) -- (-1, 2.2) -- cycle {};

    \draw
        (0, 0.2) edge[line width=2] (-1, 2.2)
        (-1, 2.2) edge[line width=2] (1, -0.8)
        (-1, 1.2) edge[line width=2] (-1, 2.2)
        (-1, 1.2) edge[line width=2] (1, -0.8);

    \filldraw[black] (-1, 2.2) circle (1pt);
    \filldraw[black] (-1, 1.2) circle (1pt);
    \filldraw[black] (0, 0.2) circle (1pt);
    \filldraw[black] (1, -0.8) circle (1pt);

    \draw
        (-1.4, 1.5) edge[line width=1] (1.4,-1.3);

    \node[below] at (0, -1.1) {$\mathcal T$};
    \node[left] at (-1.1, 0) {$\pi^* \mathcal T$};
    \node[left] at (-1.1, 1.7) {$\cone(z, (\pi^* \mathcal T)|_H)$};
    \node[right] at (1.2, -1) {$H$};
\end{tikzpicture}
    \caption{$\widecheck{P}_2^{(n + 1)} = P^{\leq H} \cup P^{\geq H}$ with its regular subdivision $\mathcal S$ glued from $\pi^* \mathcal T$ and $\cone(z, (\pi^* \mathcal T)|_H)$ along $H$.}
    \label{triangles/figure/pullbackcone}
\end{figure}

To conclude our construction, we can refine $\mathcal S$ to a unimodular triangulation $\mathcal T^{(n + 1)}$ by using the \hyperlink{triangles/pulling}{pulling refinement}.
To ensure that the triangulation is \hyperlink{background/starshaped}{star-shaped}, we pull the origin first.
See \cref{triangles/figure/pulled}.
More precisely, we can show the following.
\begin{lemma}
\label{triangles/starshapedbelow}
    Any refinement of the subdivision $\mathrm{pull}_0\parens{\pi^* \mathcal T^{(n)}}$ of $P^{\leq H}$ is \hyperlink{background/starshaped}{star-shaped}.
\end{lemma}
\begin{proof}
    The cells of $\pi^* \mathcal T^{(n)}$ are pullbacks $\pi^{*}(\sigma)$ of cells $\sigma \in \mathcal T^{(n)}$.
    Each such $\sigma$ is a unimodular simplex containing the origin (in $\R^n$) as a vertex.
    It suffices to focus our attention on any such cell $\pi^{*}(\sigma)$.
    
    The facets of the column $\pi^{*}(\sigma)$ are the ``bottom cap" $\sigma \times \{-1\}$, the ``top cap" $\pi^{*}(\sigma) \cap H$, and the pullbacks $\pi^{*}(F)$ of facets $F$ of $\sigma$.
    As $\sigma$ is a simplex, there is precisely one facet $F_0 \leq \sigma$ which does not contain the origin (in $\R^n$).
    The pullback $\pi^{*}(F_0)$ and the bottom cap $\sigma \times \{-1\}$ are the only facets of $\pi^{*}(\sigma)$ not containing the origin (in $\R^{n+1}$).
    The effect on $\pi^{*}(\sigma)$ of pulling the origin is thus to replace the facets $\sigma \times \{-1\}$ and $\pi^{*}(F_0)$ with their cones to the origin $\conv(0, \sigma \times \{-1\})$ and $\conv(0, \pi^{*}(F_0))$.
    
    $\conv(0, \sigma \times \{-1\})$ is a unimodular simplex, so it cannot be subdivided further, and it already contains the origin.
    We are thus left to show that for any subdivision of $\conv(0, \pi^{*}(F_0))$, the origin is a vertex of every cell.
    Crucially, we have that for a lattice point $m$ of $\conv(0, \pi^{*}(F_0))$, either $m = 0$ or $m$ lies in the facet $\pi^{*}(F_0)$.
    Indeed, as $\sigma$ is unimodular, its only lattice points are its vertices, which are either the origin itself or are vertices of $F_0$.
    Suppose $m$ is not in $\pi^{*}(F_0)$.
    Then $m$ must lie above $(0, \dots, 0, -1)$.
    Only $(0, \dots, 0, -1)$ and $(0, \dots, 0, 0)$ lie above this point due to \cref{triangles/highestpoints}, so $m$ must equal one of these two points.
    A representation of $m$ as a convex combination of $0$ and the vertices of $\pi^{*}(F_0)$ must contain a positive coefficient for $0$, as $m$ is not in $\pi^{*}(F_0)$.
    Thus, the last coordinate of $m$ must be strictly greater than $-1$, so $m$ cannot be $(0, \dots, 0, -1)$.
    We conclude that $m = 0$.
    
    Now, let $\Sigma$ be a cell of a subdivision of $\conv(0, \pi^{*}(F_0))$.
    For dimension reasons, $\Sigma$ cannot be wholly contained in $\pi^{*}(F_0)$.
    Hence, it must have some vertex outside of $\pi^{*}(F_0)$.
    As discussed above, the only possibility is the origin.
    We have thus shown that every cell of $\mathrm{pull}_0(\pi^* \mathcal T^{(n)})$ contains the origin as a vertex, as desired.
\end{proof}

Now, we pull all the remaining vertices in any order to define our triangulation $\mathcal T^{(n + 1)}$.
By \cref{triangles/pullingregular}, this is regular.
By \cref{triangles/pullingunimodular}, $\mathcal T^{(n + 1)}|_{P^{\leq H}}$ is unimodular.
And as shown above, $\mathcal T^{(n + 1)}|_{P^{\geq H}} = \mathcal S|_{P^{\geq H}}$ is unimodular.
Finally, by \cref{triangles/starshapedbelow}, $P^{\leq H}$ is star-shaped, and $P^{\geq H}$ is star-shaped as it is a cone of $(\pi^* \mathcal T^{(n)})|_H$, which is isomorphic to the star-shaped triangulation $\mathcal T^{(n)}$.
We have therefore found a regular, unimodular, star-shaped triangulation $\mathcal T^{(n + 1)}$ of $\widecheck{P}_2^{(n + 1)}$ that extends $\mathcal T^{(n)}$ of $\widecheck{P}_2^{(n)}$.
See \cref{triangles/figure/pulled}.
    \begin{figure}
        \centering
        \begin{minipage}{0.49\textwidth}
            \centering
	    \usetikzlibrary{quotes}

\begin{tikzpicture}[
        scale=2,
    ]
    \fill[color=red, opacity=0.6] (-1, -1) -- (0, -1) -- (0, 0) -- (-1, 1) -- cycle {};
    \fill[color=blue, opacity=0.6] (0, -1) -- (0, 0) -- (1, -1) -- cycle {};
    \draw
        (-1, -1) edge[line width=2] (1, -1)
        (-1, -1) edge[line width=2] (-1, 1)
        (0, -1) edge[line width=2] (0, 0)
        (-1, 1) edge[line width=2] (1, -1);

    \filldraw[black] (-1, -1) circle (1pt);
    \filldraw[black] (-1, 0) circle (1pt);
    \filldraw[black] (-1, 1) circle (1pt);
    \filldraw[black] (0, -1) circle (1pt);
    \filldraw[black] (0, 0) circle (1pt);
    \filldraw[black] (1, -1) circle (1pt);

    \fill[color=darkgreen, opacity=0.42] (-1, 1) -- (0, 0) -- (-1, 2) -- cycle {};
    \fill[color=pinkish, opacity=0.6] (1, -1) -- (0, 0) -- (-1, 2) -- cycle {};

    \draw
        (0, 0) edge[line width=2] (-1, 2)
        (-1, 2) edge[line width=2] (1, -1)
        (-1, 1) edge[line width=2] (-1, 2)
        (-1, 1) edge[line width=2] (1, -1);

    \filldraw[black] (-1, 2) circle (1pt);

    \draw
        (0, 0) edge[line width = 2, dashed] (-1, -1);

    \node[left] at (-0.1, 0) {$(0, 0)$};
\end{tikzpicture}
	    \caption*{(a) Pulling the origin.}
        \end{minipage}
        \begin{minipage}{0.49\textwidth}
            \centering

\begin{tikzpicture}[
        scale=2,
    ]
    \fill[color=red, opacity=0.6] (-1, -1) -- (0, -1) -- (0, 0) -- (-1, 1) -- cycle {};
    \fill[color=blue, opacity=0.6] (0, -1) -- (0, 0) -- (1, -1) -- cycle {};
    \draw
        (-1, -1) edge[line width=2] (1, -1)
        (-1, -1) edge[line width=2] (-1, 1)
        (0, -1) edge[line width=2] (0, 0)
        (-1, 1) edge[line width=2] (1, -1);

    \filldraw[black] (-1, -1) circle (1pt);
    \filldraw[black] (-1, 0) circle (1pt);
    \filldraw[black] (-1, 1) circle (1pt);
    \filldraw[black] (0, -1) circle (1pt);
    \filldraw[black] (0, 0) circle (1pt);
    \filldraw[black] (1, -1) circle (1pt);

    \fill[color=darkgreen, opacity=0.42] (-1, 1) -- (0, 0) -- (-1, 2) -- cycle {};
    \fill[color=pinkish, opacity=0.6] (1, -1) -- (0, 0) -- (-1, 2) -- cycle {};

    \draw
        (0, 0) edge[line width=2] (-1, 2)
        (-1, 2) edge[line width=2] (1, -1)
        (-1, 1) edge[line width=2] (-1, 2)
        (0, 0) edge[line width = 2] (-1, -1)
        (-1, 1) edge[line width=2] (1, -1);

    \filldraw[black] (-1, 2) circle (1pt);

    \draw
        (-1, 0) edge[line width = 2, dashed] (0, 0);

    \node[left] at (-1, 0) {$(-1, 0)$};
\end{tikzpicture}
	    \caption*{(b) Pulling the vertex $(-1, 0)$.}
        \end{minipage}

	\caption{The result of pulling the subdivision $\mathcal S$ as in \cref{triangles/figure/pullbackcone} at the origin (a) then at the remaining vertices (b). The new edges are depicted with dotted lines.}
        \label{triangles/figure/pulled}
    \end{figure}
The base case is trivial, as $\widecheck{P}_2^{(1)} = [-1,1]$.

The above construction proves \cref{triangles/triangulating2dual}, which as discussed proves \cref{triangles/triangulating2}.
Finally, as shown in \cref{section/triangles/reduction}, this proves \cref{triangles/triangulating1}.
As discussed in \cref{section/toric}, we conclude with a family of toric, projective, crepant resolutions of the weighted projective spaces $\P_1^{(n)}$ (as well as a compatible family of toric, projective, crepant resolutions of the $\P_2^{(n)}$).

	\section{Resolving the hypersurface within}
        \label{section/hypersurface}
        Consider the quasi-smooth weighted projective hypersurface
\[
    X_1^{(n)} = \{x_0^{s_0} + \dots + x_{n - 1}^{s_{n - 1}} + x_n^{d - 1} x_{n + 1} + x_{n + 1}^{d} = 0\} \subseteq \P_1^{(n + 1)},
\]
where $d = d_1^{(n + 1)} = 2s_{n} - 2$.
(Recall that a subvariety of a weighted projective space is called \emph{quasi-smooth} if its affine cone is smooth away from the origin.)
For simplicity, we write $X = X_1^{(n)}$ and $\P = \P_1^{(n + 1)}$.
$X$ is a Calabi-Yau variety of dimension $n$ from which Esser, Totaro, and Wang in \cite{largeindex} define a klt Calabi-Yau pair of dimension $n$ and a terminal Calabi-Yau variety of dimension $n + 1$ with conjecturally maximal indices.

We recall these constructions.
Let $m = (s_n - 1)(2 s_n - 3)$.
Then $\mu_m$ acts on $X$ via
\[
    \zeta [x_0, \dots, x_{n + 1}] = \left[\zeta^{d/(2s_0)} x_0 : \dots : \zeta^{d/(2s_{n-1})} x_{n - 1} : x_n : \zeta^{d/2} x_{n + 1}\right].
\]
The quotient $X/\mu_m$ is a klt Calabi-Yau pair of dimension $n$ with standard coefficients and index $m$ by \cite[Proposition 3.7]{largeindex}.
This is the conjecturally maximal index for such an $n$-fold.

For the terminal example, consider a $\mu_m$-equivariant terminalization $Z \longrightarrow X$.
Then let $E$ be an elliptic curve and $p \in E$ be a point of order $m$.
Let $\mu_m$ act on $Z \times E$ via
\[
    \zeta(x, y) = (\zeta x, y + p).
\]
Then $(Z \times E)/\mu_m$ is a terminal Calabi-Yau variety of dimension $n + 1$ with index $m$, which is conjecturally maximal among such varieties.

We seek a smooth, projective Calabi-Yau variety of dimension $n + 1$ and index $m$.
To do so, it suffices to find a $\mu_m$-equivariant, crepant, projective resolution $\widetilde{X} \longrightarrow X$ and repeat the above construction with the elliptic curve.
Furthermore, this crepant resolution $\widetilde{X}$ will be our smooth, projective Calabi-Yau variety with large sums of Betti numbers.

In \cref{section/triangles}, we constructed a toric, projective, crepant resolution $\pi: \widetilde{\P} \longrightarrow \P$.
Now, as $X$ is a quasi-smooth hypersurface in $\P$, its only singularities are the singularities of the ambient space, so we are led to believe that by resolving the singularities of the ambient space, we have resolved the singularities of $X$.
More precisely,
\begin{theorem}
\label{hypersurface/largeindex}
    Let $L$ be the linear system in $\P$ generated by the monomials $$x_0^{s_0}, \dots, x_{n - 1}^{s_{n - 1}}, x_n^{d - 1} x_{n + 1}, x_{n + 1}^{d}.$$
    Then for a generic element $X \in L$, $\pi^{-1}X \longrightarrow X$ is a projective, crepant, $\mu_m$-equivariant resolution of singularities.
\end{theorem}

To work with the canonical divisors of these hypersurfaces, we will need to use the adjunction formula.
As $\widetilde{\P}$ is smooth, the adjunction formula is known.
However, the adjunction formula is not true for every weighted projective hypersurface.
The technical condition that assures this holds is \emph{well-formedness}.
We refer to \cite{fletcher} for details.

\begin{definition}
\hypertarget{hypersurface/wellformed}
    Let $\P = \P(a_0, \dots, a_n)$ be a weighted projective space and $X \subseteq \P$ a hypersurface.
    \begin{enumerate}[label = (\roman*)]
        \item $\P$ is called \emph{well-formed} if for all $i$ we have $\gcd(a_0, \dots, \widehat{a_i}, \dots, a_n) = 1$.
        \item $X \subseteq \P$ is called \emph{well-formed} if $\P(a_0, \dots, a_n)$ is well-formed and $X$ contains no codimension $2$ singular coordinate stratum of $\P(a_0, \dots, a_n)$.
    \end{enumerate}
\end{definition}

We first note that the weights of a weighted projective space are not uniquely determined, and in fact every weighted projective space is isomorphic to a well-formed one.
For instance, using the Veronese map, we have that $\P(d a_0, \dots, d a_n) \cong \P(a_0, \dots, a_n)$.
Hence, we can reduce the weights of a weighted projective space to have no common factors.
Furthermore, if $a_0, \dots, a_n$ have no common factor and $d = \gcd(a_1, \dots, a_n)$ then $\P(a_0, \dots, a_n) \cong \P(a_0, a_1/d, \dots, a_n/d)$.
Repeatedly applying this allows the weights of any weighted projective space to be reduced to be well-formed.

Finally, as discussed in \cite[6.14]{fletcher}, the adjunction formula holds for well-formed, quasi-smooth weighted projective hypersurfaces.
We may now proceed to prove \cref{hypersurface/largeindex}.

\begin{proof}[Proof of \cref{hypersurface/largeindex}]
    Projectivity of $\pi^{-1} X \longrightarrow X$ follows from projectivity of $\pi$ and $X$.

    Let $\pi^{-1}L$ be the total transform of the linear system $L$.
    Let $\widetilde{p}$ be a base point of $\pi^{-1}L$.
    Then $p = \pi(\widetilde{p})$ is a base point of $L$, hence $p = [0 : \dots : 0 : 1 : 0]$.
    This is a smooth point of $\P$, corresponding to the cone in $\R^{n + 1}$ spanned by the unimodular simplex $\conv(e_0, \dots, e_{n - 1}, w_1)$.
    The triangulation of $P_1^{(n + 1)}$ defining $\pi$ leaves the cone $\conv(e_0, \dots, e_{n - 1}, w_1)$ unchanged, so $\pi$ is an isomorphism near $p$.
    Hence, $\widetilde{p}$ is the only base point of $\pi^{-1} L$.
    By Bertini's theorem, a generic element of $\pi^{-1} L$ is therefore smooth away from $\widetilde{p}$.
    Furthermore, $p$ is a smooth point of $X$, as $X$ is a quasi-smooth hypersurface in $\P$ and $p \in \P$ is a smooth point.
    As $\pi$ is an isomorphism near $p$, it follows that $\widetilde{p}$ is a smooth point of $\pi^{-1} X$.
    We have shown that every point of $\pi^{-1} X$ is smooth, for generic $X \in L$.

    To show crepancy of $\pi^{-1} X \longrightarrow X$, we will need to use the adjunction formula.
    As discussed above, we will need to show a generic $X \in L$ is well-formed.

    First, we show that $\P = \P\parens{d_1^{(n + 1)}/s_0, \dots, d_1^{(n + 1)}/s_{n - 1}, 1, 1}$ is well formed.
    Indeed, there are two weights equal to $1$, so if we omit only one weight the remaining $n + 1$ terms contain a $1$, and as such their greatest common divisor is $1$.

    We now show that a generic $X \in L$ does not contain any singular codimension 2 stratum of $\P$.
    Indeed, the base locus of $L$ is, as above, just the point $p = [0 : \dots : 0 : 1 : 0]$.
    So if $\dim(\P) \geq 3$, we can ensure that there is no codimension $2$ stratum whatsoever in a generic $X$.
    When $\dim(\P) = 2$, the base locus $\{p\}$ is codimension $2$, so any $X \in L$ contains $p$.
    However, as mentioned above, $p$ is a smooth point of $\P$, so there is no \emph{singular} codimension $2$ stratum for $X$ generic.

    We therefore apply the adjunction formula and write $K_X = (K_{\P} + X)|_X$.
    Using this and crepancy of $\pi: \widetilde{\P} \longrightarrow \P$, we will show that $\pi^{-1} X \longrightarrow X$ is crepant.
    Indeed,
    \begin{align*}
        K_{\pi^{-1} X} &= (K_{\widetilde{\P}} + \pi^{-1}X)|_{\pi^{-1} X}\\
        &= \pi^{-1} (K_{\P} + X)|_{X}\\
        &= \pi^{-1} K_X.
    \end{align*}

    Finally, $\pi$ is equivariant and maps the torus of $\widetilde{\P}$ isomorphically to the torus of $\P$, so the subgroup $\mu_m$ of $T_\P = T_{\widetilde{\P}}$ fixing $X$ also fixes $\pi^{-1} X$.
\end{proof}
Repeating the above construction then yields our smooth, projective Calabi-Yau variety of large index.
\begin{corollary}
    For every $n$, there is a smooth, projective Calabi-Yau variety of dimension $n + 1$ with index $(s_n - 1)(2 s_n - 3)$.
\end{corollary}
We have now proven our first main theorem \cref{intro/largeindex}.
However our examples are not fully explicit due to the genericity in \cref{hypersurface/largeindex}.
We now show that the generic examples defined there are isomorphic to an explicitly chosen one.
Suppose we have an element $X \in L$ given by
\[
    \{a_0 x_0^{s_0} + \dots + a_{n} x_n^{d - 1} x_{n + 1} + a_{n + 1}x_{n + 1}^{d} = 0\}
\]
where all $a_i$ are nonzero.
Then consider the toric automorphism of $\P$ given by 
\[
    [x_0 : \dots : x_n] \mapsto \bracket{a_0^{(1/s_0)} x_0 : \dots : (a_n a_{n + 1}^{-1/d})^{1/(d - 1)} x_{n} : (a_{n + 1}^{1/d}) x_{n + 1}}.
\]
This restricts to an isomorphism
\[
    X \xlongrightarrow{\sim} X_1^{(n)} = \{x_0^{s_0} + \dots + x_{n - 1}^{s_{n - 1}} + x_n^{d - 1} x_{n + 1} + x_{n + 1}^{d} = 0\} \subseteq \P_1^{(n + 1)}.
\]
Furthermore, as $T_{\widetilde{\P}} = T_{\P}$ this extends to a toric automorphism of $\widetilde{\P}$ which restricts to an isomorphism of their total transforms.
\[
    \pi^{-1} X \xlongrightarrow{\sim} \pi^{-1} X_1^{(n)}.
\]
As this isomorphism is given via translation by a torus element, it preserves stabilizers and is hence $\mu_m$-equivariant.
The above construction then yields \emph{explicit} smooth, projective Calabi-Yau varieties of large index, up to the choice of the order in which we pull vertices to conclude the triangulation discussed at the end of \cref{section/triangles/triangulatingP2}.

Furthermore, as mentioned in \cref{section/intro} the crepant resolutions thus constructed have extreme topological properties.
Theorem 5.1 in \cite{largeindex} computes the \emph{orbifold} Hodge numbers of $X_1^{(n)}$.
By \cite{orbifoldcrepant}, these are known to agree with the actual Hodge numbers of a crepant resolution of $X_1^{(n)}$, which we have just constructed.
Furthermore, en route to crepantly resolving $\P_1$, we crepantly resolved $\P_2$.
In \cite{largeindex}, the authors define the quasi-smooth weighted projective hypersurfaces
\[
    X_2^{(n)} := \{x_0^{s_0} + \dots + x_n^{s_n} + x_{n + 1}^{s_{n + 1} - 1} = 0\} \subseteq \P_2^{(n + 1)}
\]
which, by an argument as in \cref{hypersurface/largeindex}, admit projective, crepant resolutions.
Let $\widetilde{X}_i^{(n)}$ denote the projective, crepant resolutions of $X_i^{(n)}$ we have thus constructed.
We have then an analogous result to Theorem 5.1 in \cite{largeindex}.
\begin{theorem}
\label{hypersurface/largebetti}
    Let $n \geq 1$.
    \begin{enumerate}[label = (\roman*)]
        \item $h^{p, q}(\widetilde{X}_i^{(n)}) = 0$ for $i = 1, 2$ unless $p = q$ or $p + q = n$.
        \item The sum of the Betti numbers of $\widetilde{X}_i^{(n)}$ is
        \[
            H = 2(s_0 - 1) \cdots (s_n - 1),
        \]
        which equals the Euler characteristic when $n$ is even.
        \item If $n$ is odd then
        \[
            h^n(\widetilde{X}_i^{(n)}) = 
            \begin{cases}
                (s_0 - 1) \cdots (s_{n - 1} - 1)(2 s_n - 4) & i = 1\\
                (s_0 - 1) \cdots (s_{n} - 1) & i = 2,
            \end{cases}
        \]
        and the Euler characteristic is
        \[
            \chi(\widetilde{X}_i^{(n)}) = 
            \begin{cases}
                -(s_0 - 1) \cdots (s_{n - 1} - 1) (2 s_n - 6) & i = 1\\
                0 & i = 2.
            \end{cases}
        \]
    \end{enumerate}
\end{theorem}
\begin{proof}
    As $\widetilde{X}_i^{(n)} \longrightarrow X_i^{(n)}$ is a crepant resolution, we have that $h^{p, q}(\widetilde{X}_i^{(n)}) = h^{p, q}(X_i^{(n)})$ by \cite{orbifoldcrepant}.
    The result then follows from Theorem 5.1 in \cite{largeindex}, where they compute these exact values as the orbifold Hodge/Betti numbers of $X_i^{(n)}$.
\end{proof}
This shows our second main theorem, \cref{intro/largebetti}.
Additionally, see \cref{hypersurface/hodgediamonds} for the Hodge diamonds of our $3$ and $4$ dimensional examples, as computed in \cite{largeindex}.

\begin{figure}
    \begin{minipage}{0.5\textwidth}
        \centering
        \begin{tabular}{c c c c c c c}
            & & & 1 & & &\\
            & & 0 & & 0 & &\\
            & 0 & & 11 & & 0 &\\
            1 & & 491 & & 491 & & 1\\
            & 0 & & 11 & & 0 &\\
            & & 0 & & 0 & &\\
            & & & 1 & & &
        \end{tabular}
        \caption*{(a) $\widetilde{X}_1^{(3)}$}
    \end{minipage}
    \begin{minipage}{0.5\textwidth}
        \centering
        \begin{tabular}{c c c c c c c}
            & & & 1 & & &\\
            & & 0 & & 0 & &\\
            & 0 & & 251 & & 0 &\\
            1 & & 251 & & 251 & & 1\\
            & 0 & & 251 & & 0 &\\
            & & 0 & & 0 & &\\
            & & & 1 & & &
        \end{tabular}
        \caption*{(b) $\widetilde{X}_2^{(3)}$}
    \end{minipage}

    \vspace{1cm}

    \resizebox{\linewidth}{!}{%
        \begin{minipage}{0.5\textwidth}
            \centering
            \begin{tabular}{c c c c c c c c c}
                & & & & 1 & & & &\\
                & & & 0 & & 0 & &  \\
                & & 0 & & 252 & & 0 & &\\
                & 0 & & 0 & & 0 & & 0\\
                1 & & 303148 & & 1213644 & & 303148 & & 1\\
                & 0 & & 0 & & 0 & & 0\\
                & & 0 & & 252 & & 0 & &\\
                & & & 0 & & 0 & && \\
                & & & & 1 & & & &
            \end{tabular}
            \caption*{(c) $\widetilde{X}_1^{(4)}$}
        \end{minipage}
        \hspace{0.5cm}
        \begin{minipage}{0.5\textwidth}
            \begin{tabular}{c c c c c c c c c}
                & & & & 1 & & & &\\
                & & & 0 & & 0 & &  \\
                & & 0 & & 151700 & & 0 & &\\
                & 0 & & 0 & & 0 & & 0\\
                1 & & 151700 & & 1213644 & & 151700 & & 1\\
                & 0 & & 0 & & 0 & & 0\\
                & & 0 & & 151700 & & 0 & &\\
                & & & 0 & & 0 & && \\
                & & & & 1 & & & &
            \end{tabular}
            \caption*{(d) $\widetilde{X}_2^{(4)}$}
        \end{minipage}
    }
    \caption{The Hodge diamonds of the crepant resolutions $\widetilde{X}_i^{(n)}$ for $i = 1, 2$ and $n = 3, 4$.}
    \label{hypersurface/hodgediamonds}
\end{figure}

Recall \cref{intro/smoothbetti}, which states that these examples have the largest Betti number sums and for odd $n$, that the Euler characteristics are as negative as possible.
There is a notable omission here.
In \cite{largeindex}, the authors also constructed a third quasi-smooth weighted projective hypersurface $X_3^{(n)} \subseteq \P_3^{(n + 1)}$, which they conjectured to have the largest possible orbifold Euler characteristic among projective varieties with quotient singularities and trivial canonical class in odd dimensions.
$X_3^{(n)}$ is mirror to $X_1^{(n)}$ (also, $X_2^{(n)}$ is self-mirror), and its orbifold Euler characteristic is negative to that of $X_1^{(n)}$.
We have not constructed a crepant resolution of $X_3^{(n)}$ so we are as of yet unable to extend these computations and conjectures to actual Betti numbers in the smooth case.

\newpage

\emergencystretch=3em
\printbibliography

@article {unimodular,
    AUTHOR = {Haase, Christian and Paffenholz, Andreas and Piechnik, Lindsay
              C. and Santos, Francisco},
     TITLE = {Existence of unimodular triangulations---positive results},
   JOURNAL = {Mem. Amer. Math. Soc.},
  FJOURNAL = {Memoirs of the American Mathematical Society},
    VOLUME = {270},
      YEAR = {2021},
    NUMBER = {1321},
     PAGES = {v+83},
      ISSN = {0065-9266,1947-6221},
      ISBN = {978-1-4704-4716-8; 978-1-4704-6530-8},
   MRCLASS = {52B20 (13F20 13P10 14M25)},
  MRNUMBER = {4277268},
MRREVIEWER = {Margaret\ M.\ Bayer},
       DOI = {10.1090/memo/1321},
       URL = {https://doi.org/10.1090/memo/1321},
}

@article{largeindex,
	title = {Calabi-{Yau} varieties of large index},
	url = {https://arxiv.org/abs/2209.04597},
	%doi = {10.48550/ARXIV.2209.04597},
%   urldate = {2024-08-29},
	author = {Esser, Louis and Totaro, Burt and Wang, Chengxi},
	year = {2022},
	note = {To appear in \emph{Alg. Geom}},
}

@book{toric,
	title = {Introduction to {Toric} {Varieties}. ({AM}-131)},
	isbn = {9780691000497},
%	urldate = {2024-09-06},
	publisher = {Princeton University Press},
	author = {Fulton, William},
	year = {1993},
}

@misc{nocrepant,
    title = {Quotient singularities with no crepant resolution?},
    author = {Sandor Kovacs},
    howpublished = {MathOverflow},
%    note = {URL:https://mathoverflow.net/q/66702 (version: 2022-12-13)},
    url = {https://mathoverflow.net/q/66702},
}

@article{lci,
	title = {All toric local complete intersection singularities admit projective crepant resolutions},
	volume = {53},
%	issn = {0040-8735, 2186-585X},
%	url = {https://projecteuclid.org/journals/tohoku-mathematical-journal/volume-53/issue-1/All-toric-local-complete-intersection-singularities-admit-projective-crepant-resolutions/10.2748/tmj/1178207533.full},
	doi = {10.2748/tmj/1178207533},
	number = {1},
%	urldate = {2023-09-19},
	journal = {Tohoku Mathematical Journal},
	author = {Dais, Dimitrios I. and Haase, Christian and Ziegler, Günter M.},
%	month = jan,
	year = {2001},
	note = {Publisher: Tohoku University, Mathematical Institute},
	pages = {95--107},
}

@article{smoothindices,
    author = {Masamura, Yuto},
    title = {Relations Between Indices of Calabi–Yau Varieties and Pairs},
    journal = {International Mathematics Research Notices},
    volume = {2025},
    number = {10},
    pages = {rnaf114},
    year = {2025},
    month = {05},
    issn = {1073-7928},
    doi = {10.1093/imrn/rnaf114},
}

@phdthesis{indexconjecture,
	title = {Complements on log canonical Fano varieties and index conjecture of log Calabi-Yau varieties},
	url = {https://files01.core.ac.uk/download/pdf/337604897.pdf},
	author = {Xu, Yanning},
%	month = may,
	year = {2020},
	school = {University of Cambridge, Trinity College},
	type = {Doctoral thesis},
}

@article{orbifold,
  	title = {A {New} {cohomology} {theory} of {orbifold}},
	volume = {248},
	copyright = {http://www.springer.com/tdm},
	%issn = {0010-3616, 1432-0916},
	%url = {http://link.springer.com/10.1007/s00220-004-1089-4},
	doi = {10.1007/s00220-004-1089-4},
	number = {1},
	journal = {Communications in Mathematical Physics},
	author = {Chen, Weimin and Ruan, Yongbin},
	year = {2004},
	pages = {1--31},
}

@article{orbifoldcrepant,
	title={Twisted jets,
	motivic measures and orbifold cohomology},
	volume={140},
	DOI={10.1112/S0010437X03000368},
	number={2},
	journal={Compositio Mathematica},
	author={Yasuda,
	Takehiko},
	year={2004},
	pages={396–422}
}

@article{kawamata86,
	title = {On the plurigenera of minimal algebraic 3-folds with $K \equiv 0$},
	volume = {275},
%	issn = {1432-1807},
	%url = {https://doi.org/10.1007/BF01459135},
	doi = {10.1007/BF01459135},
%	language = {en},
	number = {4},
%	urldate = {2024-12-10},
	journal = {Mathematische Annalen},
	author = {Kawamata, Yujiro},
%	month = dec,
	year = {1986},
	pages = {539--546},
}

@article{morrison86,
	title = {A remark on {Kawamata}'s paper “{On} the plurigenera of minimal algebraic 3-folds with $K \equiv 0$”},
	volume = {275},
%	issn = {1432-1807},
	%url = {https://doi.org/10.1007/BF01459136},
	doi = {10.1007/BF01459136},
%	language = {en},
	number = {4},
%	urldate = {2024-12-10},
	journal = {Mathematische Annalen},
	author = {Morrison, David R.},
%	month = dec,
	year = {1986},
	pages = {547--553},
}

@article{ps09,
	title = {Towards the second main theorem on complements},
	volume = {18},
	issn = {1056-3911, 1534-7486},
	url = {https://www.ams.org/jag/2009-18-01/S1056-3911-08-00498-0/},
	doi = {10.1090/S1056-3911-08-00498-0},
	number = {1},
	urldate = {2026-01-09},
	journal = {Journal of Algebraic Geometry},
	author = {Prokhorov, Yu. and Shokurov, V.},
	year = {2008},
	pages = {151--199},
}

@article{jiang21,
	title = {A gap theorem for minimal log discrepancies of non-canonical singularities in dimension three},
	volume = {30},
	doi = {10.1090/jag/759},
	number = {4},
	journal = {Journal of Algebraic Geometry},
	author = {Jiang, Chen},
	year = {2021},
	pages = {759--800},
}

@article{jl21,
	title = {Boundedness of log pluricanonical representations of log {Calabi}-{Yau} pairs in dimension 2},
	volume = {15},
	doi = {10.2140/ant.2021.15.547},
	number = {2},
	journal = {Algebra \& Number Theory},
	author = {Jiang, Chen and Liu, Haidong},
	year = {2021},
	pages = {545--567},
}

@article{sylvester,
	title = {On {Kellogg}'s {Diophantine} {problem}},
	volume = {29},
%	issn = {0002-9890},
%	url = {https://www.jstor.org/stable/2299023},
	doi = {10.2307/2299023},
	number = {10},
%	urldate = {2024-12-10},
	journal = {The American Mathematical Monthly},
	author = {Curtiss, D. R.},
	year = {1922},
%	note = {Publisher: [Taylor \& Francis, Ltd., Mathematical Association of America]},
	pages = {380--387},
}

@misc{oeis,
    Author = {{OEIS Foundation Inc.}},
    Note = {Published electronically at \url{http://oeis.org}},
    Title = {The {O}n-{L}ine {E}ncyclopedia of {I}nteger {S}equences},
}

@article{birkar2025,
      title={Singularities on Fano fibrations and beyond}, 
      author={Caucher Birkar},
      year={2025},
      eprint={2305.18770},
      archivePrefix={arXiv},
      primaryClass={math.AG},
      url={https://arxiv.org/abs/2305.18770}, 
}

@book{polytope,
	address = {New York, NY},
	series = {Graduate {Texts} in {Mathematics}},
	title = {Lectures on {Polytopes}},
	volume = {152},
	copyright = {http://www.springer.com/tdm},
	isbn = {9781461384311},
%	url = {http://link.springer.com/10.1007/978-1-4613-8431-1},
%	urldate = {2025-01-21},
	publisher = {Springer New York},
	author = {Ziegler, Günter M.},
	year = {1995},
	doi = {10.1007/978-1-4613-8431-1},
}

@article {k3,
    AUTHOR = {Machida, Natsumi and Oguiso, Keiji},
     TITLE = {On {K3} surfaces admitting finite non-symplectic group
              actions},
   JOURNAL = {J. Math. Sci. Univ. Tokyo},
  FJOURNAL = {The University of Tokyo. Journal of Mathematical Sciences},
    VOLUME = {5},
      YEAR = {1998},
    NUMBER = {2},
     PAGES = {273--297},
      ISSN = {1340-5705},
}

@incollection {fletcher,
    AUTHOR = {Iano-Fletcher, A. R.},
     TITLE = {Working with weighted complete intersections},
 BOOKTITLE = {Explicit birational geometry of 3-folds},
    SERIES = {London Math. Soc. Lecture Note Ser.},
    VOLUME = {281},
     PAGES = {101--173},
 PUBLISHER = {Cambridge Univ. Press, Cambridge},
      YEAR = {2000},
      ISBN = {0-521-63641-8},
}

@article {introwps,
    AUTHOR = {Beltrametti, Mauro and Robbiano, Lorenzo},
     TITLE = {Introduction to the theory of weighted projective spaces},
   JOURNAL = {Exposition. Math.},
  FJOURNAL = {Expositiones Mathematicae. International Journal for Pure and
              Applied Mathematics},
    VOLUME = {4},
      YEAR = {1986},
    NUMBER = {2},
     PAGES = {111--162},
      ISSN = {0723-0869},
}

@manual{sagemath,
    Key          = {SageMath},
    Author       = {{The Sage Developers}},
    Title        = {{S}ageMath, the {S}age {M}athematics {S}oftware {S}ystem ({V}ersion 10.6)},
    note         = {{\tt https://www.sagemath.org}},
    Year         = {2025},
}

\end{document}